\documentclass[12pt,twoside,english,reqno,a4paper]{amsart}
\usepackage{listings,graphicx,amsmath,varioref,amscd,amssymb,color,bm,stmaryrd,amsthm,amsfonts,graphics,geometry,latexsym,pgf,pst-all}
\theoremstyle{plain}
\usepackage{esint}
\usepackage{amsthm}

\theoremstyle{plain}
\newtheorem{theorem}{Theorem}[section]
\newtheorem{proposition}[theorem]{Proposition}

\newtheorem{lemma}[theorem]{Lemma}

\newtheorem{corollary}[theorem]{Corollary}

\usepackage{geometry}
\usepackage[colorlinks=false]{hyperref}

\theoremstyle{definition}
\newtheorem{defin}[theorem]{Definition}

\newtheorem{remark}[theorem]{Remark}

\theoremstyle{remark}

\def\bk{\color{black}}

\numberwithin{equation}{section}

\def\into{\int_{\Omega}}

\DeclareMathOperator{\sgn}{sgn}

\DeclareMathOperator{\N}{\mathbb{N}}

\newcommand{\car}[1]{\raise1pt\hbox{$\chi$}_{#1}}

\newcommand{\DM }{\mathcal{DM}^\infty }

\newcommand{\res}{\!\!\mathop{\hbox{
			\vrule height 7pt width .5pt depth 0pt
			\vrule height .5pt width 6pt depth 0pt}}
	\nolimits}
\begin{document}
	\title[Transparent media with $L^1$--data]{Stationary equation of the relativistic heat diffusion in transparent media having $L^1$--data}
	
	\author[F. Balducci]{Francesco Balducci}
	\author[S. Segura de Le\'on]{Sergio Segura de Le\'on}
	\address{Francesco Balducci
		\hfill \break\indent
		Dipartimento di Scienze di Base e Applicate per l' Ingegneria, Sapienza Universit\`a di Roma
		\hfill \break\indent
		Via Scarpa 16, 00161 Roma, Italy}
	\email{\tt francesco.balducci@uniroma1.it}
	\address{Sergio Segura de Le\'on
		\hfill \break\indent
		Departament d’An\`alisi Matem\`atica
		Universitat de Val\`encia
		\hfill \break\indent
		Dr. Moliner 50,
		46100 Burjassot Val\`encia, Spain}
	\email{\tt sergio.segura@uv.es}


	\begin{abstract}
		Our objective is to prove existence of a solution to the Dirichlet problem for an equation arising in the theory of radiation hydrodynamics to deal with the radiating energy in transparent media. We study its stationary equation with $L^1$--datum in a bounded domain. This problem was addressed in \cite{BOPS} for regular data (data belonging to $L^N(\Omega)$) and a bounded solution was obtained. In our framework, the proof of existence is far from trivial since the solution sought cannot be bounded. Consequently, the Anzellotti theory of pairings does not apply and we have to use new developments to introduce the meaning of solution. We also study the regularity of solutions when data belong to $L^p(\Omega)$, with $1<p<N$. Our result is coherent with the regularity found in \cite{BOPS}.	
		
	\end{abstract}
	
	\maketitle
	\tableofcontents
	
	\section{Introduction}
	
	This paper is devoted to study existence and regularity of an elliptic equation driven by the nonlinear transparent media operator with $L^1$--datum in a bounded domain. More precisely, we study the Dirichlet problem for
	\begin{equation}\label{es}
	\begin{cases}
	\displaystyle	- \operatorname{div}\left(\displaystyle u^m \frac{Du}{|Du|}\right)=f & \text{in $\Omega$,}\\ u=0 & \text{on $\partial \Omega$,}
	\end{cases}	
	\end{equation}
	where $m>0$ and $f$ is a nonnegative datum belonging to $L^1(\Omega)$.
	
	The evolution equation
	\begin{equation}\label{ev-1}
		u_t=c \operatorname{div}\left(\displaystyle u \frac{Du}{|Du|}\right),
	\end{equation}
	was proposed by J.R. Wilson using a flux limiter to deal with the radiation in transparent media with constant speed of propagation $c$, which is the highest admissible speed for transport of radiation (see \cite[Chapter IV]{MM}). The Cauchy problem for \eqref{ev-1} was studied in \cite{ACMM} where it is called the relativistic heat equation in transparent media. This analysis was extended to a class of flux limited diffusion equations including 
	\begin{equation}\label{ev-2}
		u_t=c \operatorname{div}\left(\displaystyle u^m \frac{Du}{|Du|}\right)\,,
	\end{equation}
	with $m\ge1$, in \cite{C} (we refer to \cite{CCCSS} for a survey of this class of diffusion equations). Further developments for \eqref{ev-2}, with $m\ne0$, can be found in \cite{GMP, GMP2} where is designated as the relativistic porous medium equation. This name derives from being $\operatorname{div}\left(\displaystyle u^m \frac{Du}{|Du|}\right)$ the formal limit of $\Delta_p\big(u^{m/(p-1)}\big)$ when $p$ goes to 1. Here, as usual, $\Delta_p u=\operatorname{div}\left(\displaystyle |\nabla u|^{p-2}\nabla u\right)$ stands for the $p$-Laplacian operator. The limiting case, $m=0$, is the Total Variation Flow
	\begin{equation*} 
		u_t= \operatorname{div}\left(\displaystyle \frac{Du}{|Du|}\right)\,,
	\end{equation*}
	which has been analyzed by a large number of authors due to its application to image processing (see \cite{ROF,AFCM,BCRS,Meyer,Sapiro}).
	
	The stationary equation of the Total Variation Flow is the elliptic equation governed by the 1-Laplacian. It satisfies some special features. For instance, in Lebesgue spaces a solution can only be found when $f\in L^N(\Omega)$ and its norm is small enough (see \cite{CT} and \cite{MST0}). Alternatively, unless $f\in L^N(\Omega)$ verifies a size condition, the solution found as a limit of solutions to $p$-Laplacian problems, by letting $p$ go to 1, is not a.e. finite (see \cite[Theorem 4.2]{MST}). Moreover, in general, solutions present a jump part.

	Observing that the nonlinear transparent media operator can formally be expressed through the 1-Laplacian, equation \eqref{es} is formally stated as
	\begin{equation}\label{formal}
		-\Delta_1 u=m\frac{|Du|}{u}+\frac{f}{u^m}.
	\end{equation}
	From this way of writing it, some features can be expected. The presence of the lower order term of order one should imply that the solution has not jump part. Indeed, equations involving the 1-Laplacian having a gradient term with natural growth have been widely considered (see \cite{ADS, MS, LS, LOPS, DGOP, GOP, BOP, B1}). A common property to all of them is that the solutions do not have a jump part. On the other hand, the term $\frac{f}{u^m}$ has a regularizing effect which leads to obtaining a solution without the need to impose a smallness condition (see \cite{DGS, DGOP}). Hence, these results contrast with the features of the limiting case $m=0$, the equation driven by the $1$-Laplacian.
	
A first study of the elliptic nonlinear transparent media equation in a bounded domain can be found in \cite{BOPS}. Summarizing, the main results found in this paper are:
\begin{itemize}
	\item  The solution does not have a jump part and is always bounded, independently of the size of the norm of $f$. This is contrary to what happens in the context of the $1$-Laplacian, where a condition on the size of the norm is necessary.
	\item  Moreover, if $f \not\equiv 0$ then $u \not\equiv 0$, contrary to the $1$-Laplacian case where trivial solutions can exist even in the presence of a non-trivial datum $f$. (For further details on the $1$-Laplacian, see \cite{CT,MST0}).
\end{itemize}
	Therefore, its solutions shares the main properties than solutions to the formal equation \eqref{formal}. Nevertheless, these properties are much more difficult to demonstrate in \cite{BOPS}. Let us briefly explain the reasons. The energy space to handle problems involving operators having a growth similar to the 1-Laplacian or the nonlinear  transparent media operator is the space of functions of bounded variation. It consists of $L^1$ functions whose distributional gradient is a Radon measure.  The concept of solution to the 1-Laplacian was introduced in \cite{D, ABCM}. The meaning of the quotient $\frac{Du}{|Du|}$ (quotient between two Radon measures) is achieved through a vector field $z\in L^\infty(\Omega)^N$ which satisfies $\| z\|_{L^\infty(\Omega)^N}\le1$ and $(z, Du)=|Du|$. The pairing $(z, Du)$ is a generalization of the dot product; it was introduced by Anzellotti (see \cite{A}) to define the product of a bounded vector field whose divergence is a function (or even a Radon measure) and the gradient of a $BV$-function, and obtain a Green's formula. This setting hardly changes when considering equations with a gradient term. However, the nonlinear  transparent media operator brings more difficulties. Now two vector fields have to be regarded. One of them, say $w$, plays the role of the quotient $\frac{Du}{|Du|}$ while the other, say $z$, also involves the solution: $z=u^mw$. It holds $\|w\|_{L^\infty(\Omega)^N}\le 1$ and the identity $(w, Du)=|Du|$ should be fulfilled, but that pairing has no meaning in Anzellotti's theory since we do not know that $\operatorname{div}w$ is a Radon measure. That is why more work is needed to obtain results similar to those of the 1-Laplacian with a gradient term.
	
	In our $L^1$ framework these difficulties are even more awkward. Indeed, solutions cannot be expected to be bounded, so that Anzellotti's theory is not available. Instead, a new theory which deals with pairings of unbounded vector fields must be applied. This theory has already been introduced in \cite{CDCS}, but it is not as simple as that by Anzellotti. Indeed, more preliminaries are required to show this theory in our context and more troubles to carry out our arguments.
	
	In this context, the main difficulty arises from the fact that the normal component of our vector field $z$ on a smooth hypersurface possesses a notion of trace only as a distribution and not as a function (for more details, see Remark \ref{abs-cont-Haus}).
	
	As shown in \cite[Lemma 5.3]{BOPS}, the notion of a weak trace for the normal component of $z$ on a hypersurface is used to prove that the jump part of $u$ is empty, that is, $u \in DTBV(\Omega)$. Therefore, in our more general case, we proceed by truncating the field using the characteristic function $\chi_{\{u \le k\}}$, with $k>0$. Thus, in Lemma \ref{lemma delle misure} is demonstrated that for any $z \chi_{\{ u \le k\}} \in \mathcal{DM}^\infty(\Omega)$, there is a trace for the normal component of this field on regular hypersurfaces and it is a bounded function. This, in turn, allows us to prove that $u \in DTBV(\Omega)$ even in this generalized context (for further details, see the proof of Lemma \ref{lem u j v}).
	
	Regarding the trace of the normal component of $z$ on the boundary, we will show in Lemma \ref{lem trac z in l1} that $z$ admits a trace in $L^1(\partial \Omega)$ and it can be obtained  as the limit of the traces $[z\chi_{\{u \le k\}},\nu]$ as $k\to \infty$.
	
	Finally, we highlight an important difference with the problem
	$$\begin{cases}
		\displaystyle -\Delta_1 u =f & \text{in $\Omega$,}\\
		u=0 & \text{on $\partial \Omega$,}
	\end{cases}$$
	with $0 \le f \in L^1(\Omega)$ studied in \cite{MST}. For problem \eqref{es}, we prove that the solution $u$ is finite almost everywhere in $\Omega$, in contrast to what happens in the context of the $1$-Laplacian operator.

	This paper is organized as follows.
	\begin{itemize}
		\item Section \ref{sec2}: We provide our notation, the functional spaces that we will use for the proofs and some useful tools.
		\item Section \ref{sec3}:  We present the Anzellotti pairing theory and its generalization to the case of unbounded vector fields.
		\item Section \ref{sec4}:  We introduce the notion of solution and state the main result.
		\item Section \ref{sec5}:  We show the proof of the main result.
		\item Section \ref{sec6}:  We provide a regularity result.
	
	\end{itemize}
	\black
	
	\section{Preliminaries}\label{sec2}
	\subsection{Fundamental Notation}
	In this paper,  $\mathcal L^N$ represents the $N$-dimensional Lebesgue measure, while $\mathcal{H}^{N-1}$ denotes the $(N-1)$--dimensional Haussdorff measure. Nevertheless, we will usually write $|E|$ instead of $\mathcal L^N(E)$ for the Lebesgue measure of a set $E$ and $dx$ instead of $d\mathcal L^N$.
	\\
	Henceforth, $\Omega \subset \mathbb{R}^N$, with $N \geq 2$, stands for a bounded open set with a boundary which is Lipschitz continuous. Thus, an outward normal unit vector
	$\nu(x)$ is defined for $\mathcal H^{N-1}$-a.e. $x \in\partial\Omega$.
	
	$\mathcal{B}(\Omega)$ signifies the space of Borel measurable functions defined on $\Omega$ where a.e. equal functions are identified, instead, $\mathcal{B}(\Omega; [0,1])$ denotes the space of functions $v \in \mathcal{B}(\Omega)$ taking values in the interval $[0,1]$. This notation will be adopted henceforth to make explicit the range in function spaces.
	\\We denote by $\mathcal{M}(\Omega)$ the space of Radon measure with finite total variation over $\Omega$.
	\medskip
	\\We need to consider certain truncation functions that are useful in what follows. Let $-\infty \le a < b \le \infty$, and let $T_{a,b} : \mathbb{R} \to \mathbb{R}$ be given by
	\begin{equation*}
		T_{a,b}(s):=\begin{cases}
			a & \text{if $s<a$,}\\
			s & \text{if $a\le s \le b$,}\\
			b & \text{if $s < b$.}
		\end{cases}
	\end{equation*}
	Sometimes, we need to consider the limiting function $T_{a,\infty}$, which has an obvious meaning.
			\\In the special case where $b = -a$, we obtain the standard truncation function $T_b : \mathbb{R} \to \mathbb{R}$, defined as
	\begin{equation}\label{T_k}
		T_b(s):=\begin{cases}
			s & \text{if $|s|\le b$,} \\
			\sgn(s)b & \text{if $|s|>b$.}
		\end{cases}
	\end{equation}
	Furthermore, for our purposes, it is essential to introduce the following function, let $j, \varepsilon>0$, we define $h_{j,\varepsilon}: \mathbb{R} \to \mathbb{R}$ by
	\begin{equation}\label{hkeps}
		h_{j, \varepsilon}(s):=\begin{cases}
			0 & \text{if $s \ge j+\varepsilon,$}\\
			1 & \text{if $s \le j$,}\\
			\frac{j+\varepsilon-s}{\varepsilon} & \text{if $j<s<j+\varepsilon$.}
	\end{cases} \end{equation}
	\medskip
	\\To maintain clarity and avoid ambiguity, we will frequently adopt the following convention:
	$$\int_{\Omega} f := \int_{\Omega} f(x) \, \ensuremath{d}x,$$
	the Hausdorff measure will always be made explicit
	and, if $\mu$ is another Radon measure, then we will write
	\begin{equation*}
		\int_\Omega f \mu := \int_\Omega f \, \ensuremath{\mathrm d}\mu.
	\end{equation*}
	\medskip
	\\Unless explicitly stated otherwise, we use $\mathcal{C}$ to represent various positive constants, the values of which may vary from one line to the next and, occasionally, even within the same line. These constants depend only on the problem's data and are independent of the indices of the sequences introduced throughout. Additionally, for simplicity, we will not relabel an extracted convergence subsequence when no risk of confusion arises.
	
	The following result is a well-known consequence of Egorov's Theorem. We state it for further reference.
	\begin{lemma}\label{lemmino di Sergio}
		Let us consider $(X, \mu)$ a measurable space with $\mu$ a nonnegative finite measure. If  $f_n \in L^1(X,\mu)$ defines a sequence such that $f_n \rightharpoonup f$ weakly in $L^1(X,\mu)$ and $g_n \in L^1(X, \mu)$ defines another one such that $g_n \to g$ $\mu$-a.e. in $X$ and it is uniformly bounded, that is, $$|g_n| \le \mathcal{C},$$ where $\mathcal{C}>0$ is a constant independent of $n \in \mathbb{N}$, then it holds $$\int_X fg \, \mu = \lim_{n \to \infty } \int_X f_n g_n \, \mu.$$
	\end{lemma}
	
	\subsection{Marcinkiewicz space}
	In this subsection. we introduce the Marcinkiewicz spaces, Let us consider $q > 0$, and define $L^{q,\infty}(\Omega)$ as the space of measurable functions $v: \Omega \to \mathbb{R}$ for which there exists a constant $\mathcal{C} > 0$ such that
	$$\left|\left\{|v|>t\right\} \right| \le \frac{\mathcal{C}}{t^q}\quad \text{for all }t>0.$$
	So the following quantity
	$$[v]_{L^{q,\infty}(\Omega)} := \inf \left\{ \mathcal{C} > 0 \, : \, \left| \left\{ |v| > t \right\} \right| \le \frac{\mathcal{C}}{t^q} \right\},$$
	is always finite. We remark that, for $1 \le q \le \infty $, this does not give rise to a genuine norm on $L^{q,\infty}(\Omega)$; rather, it is just a quasi-norm.
	We recall that the following inclusions hold:
	$$L^p(\Omega) \subset L^{q, \infty}(\Omega) \subset L^r(\Omega),$$
	for every $1 \leq r < q \le p\le \infty$. For further information, we refer to \cite[Appendix]{DL}.
	\subsection{Essential notions on $BV$ and $TBV$ functions}\label{subs BV}
	We commence by presenting the notion of functions of bounded variation. For further details, we direct the reader to \cite[Chapter 3]{AFP} or \cite[Chapter 10]{ABM}.
	\\We define the space
	$$BV(\Omega) := \left\{ v \in L^1(\Omega) \, : \, Dv \in \mathcal{M}(\Omega)^N \right\},$$
	where $Dv$ denotes the distributional gradient. The total variation of this vector measure is represented by  $|Dv|$ and defined as
	$$\int_\Omega|Dv| := \sup \left\{ \int_\Omega v \operatorname{div} \psi \, : \, \psi \in C^1_c(\Omega)^N, \, \|\psi\|_{L^\infty(\Omega)^N} \le 1 \right\}.$$
	It is proven that functions of bounded variation possess a trace on the boundary $v$ which belongs to $L^1(\partial \Omega)$. The space $BV(\Omega)$ can be endowed with the norm $$\|v\|_{BV(\Omega)}:=\int_{\partial \Omega}v\, d\mathcal H^{N-1}+\int_\Omega |Dv|\,,$$
	with which is a Banach space.
	\medskip
	\\We now recall the concept of approximate limit. Given $v \in BV(\Omega)$, we say that $\tilde{v}(x) \in \mathbb{R}$ is the approximate limit of $v$ at $x$ if
	$$\lim_{\rho \to 0} \frac{1}{|B_\rho(x)|} \int_{B_\rho(x)} |v(y) - \tilde{v}(x)| \, \ensuremath{d}y = 0.$$
	We denote by $S_v$ the set of  points where an approximate limit of the function $v$ does not exist.
	\\We say that $x \in \Omega$ is a jump point for the function $v$ if there exists a unit vector $\nu(x)$ and two distinct values $v^e(x), v^i(x) \in \mathbb{R}$ such that
	$$\lim_{\rho \to 0} \frac{1}{|B_\rho^e(x)|}\int_{B_\rho^e(x)} |v(y)-v^e(x)|\, dy=0, $$
	$$\lim_{\rho \to 0} \frac{1}{|B_\rho^i(x)|}\int_{B_\rho^i(x)} |v(y)-v^i(x)|\, dy=0,$$
	where
	$$B_\rho^e(x):= \left\{y \in B_\rho(x) \,:\, (y-x)\cdot \nu(x) >0\right\},$$
	$$B_\rho^i(x):= \left\{y \in B_\rho(x) \,:\, (y-x)\cdot \nu(x) <0\right\}.$$
	The symbol $J_v$ will be used to denote the set of jump points, which satisfies $J_v\subset S_v$. Furthermore, if $v \in BV(\Omega)$, it can be shown that $\mathcal{H}^{N-1}(S_v \setminus J_v) = 0$ (see \cite[Theorem 3.78]{AFP}).
	
	Moreover, $J_v$ is an $\mathcal{H}^{N-1}$-rectifiable set, and it is possible to define an orientation $\nu_v(x)$ for $\mathcal{H}^{N-1}$-a.e. $x \in J_v$.
	\medskip
	\\Given $v \in L^1(\Omega)$, we define its precise representative as the function $v^* : \Omega \setminus (S_v \setminus J_v) \to \mathbb{R}$ given by
	$$v^*(x):=\begin{cases}
		\tilde{v}(x) & \text{if $x \in \Omega\backslash S_v$,}\\
		\frac{v^e(x)+v^i(x)}{2} & \text{if $x \in J_v$,}
	\end{cases}$$
	moreover, we emphasize that if $v \in BV(\Omega)$, then $v^*$ is determined $\mathcal{H}^{N-1}$-a.e., since $\mathcal{H}^{N-1}(S_v \setminus J_v) = 0$.
	\medskip
	\\For a function $v \in BV(\Omega)$, by the Radon-Nikodým Theorem, we can decompose $Dv$ into an absolutely continuous part with respect to the Lebesgue measure, denoted $D^a v$, and a singular part, denoted $D^s v$. The singular part can further be decomposed into the jump derivative $D^j v=D^s v\res{J_v}$ and the Cantor derivative $D^c v=D^s v\res{\Omega\backslash S_v}$. Summarizing, we have
	$$ Dv = D^a v + D^s v = D^a v + D^j v + D^c v.$$
	We point out that, when $D^j v = 0$, or equivalently, $\mathcal{H}^{N-1}\left(J_v\right)=0$, we will write $v$ instead of $v^*$ because there is no ambiguity. We will denote by $DBV(\Omega)$ the space of functions $v \in BV(\Omega)$ such that $D^j v = 0$.
	\medskip
	\\We now recall a result on lower semi-continuity for functions of bounded variation (see \cite[Proposition 3.6]{AFP}).
	\begin{lemma} 
		Let us consider a sequence $v_n \in BV(\Omega)$ such that $v_n \to v$ strongly in $L^1(\Omega)$ with $v \in BV(\Omega)$. Then
		\begin{equation}\label{sci grad con bordo}
			\int_{\partial \Omega} v \varphi \, \ensuremath{\mathrm d}\mathcal{H}^{N-1} + \into |Dv| \varphi \le \liminf_{n \to \infty} \int_{\partial \Omega} v_n \varphi \, \ensuremath{\mathrm d}\mathcal{H}^{N-1} + \into |Dv_n|\varphi,
		\end{equation}
		for all $0 \le \varphi \in C^1(\overline{\Omega})$.
	\end{lemma}
	As a straightforward consequence, if $0\le \varphi \in C^1_c(\Omega)$, then
	\begin{equation}\label{sci grad}
		\int_\Omega |Dv|\varphi \le \liminf_{n \to \infty} \int_\Omega |Dv_n| \varphi \,.
	\end{equation}
	
	We now present Sobolev embedding Theorem (\cite[Proposition 3.23]{AFP}).
	\begin{theorem}\label{embedding}
		The embedding $BV(\Omega) \hookrightarrow L^{\frac{N}{N-1}}(\Omega)$ is continuous, i.e. there exists a constant $\mathcal{S}_1>0$ such that
		\begin{equation}\label{sob emb}
			\|v\|_{L^{\frac{N}{N-1}}(\Omega)} \le \mathcal{S}_1\|v\|_{BV(\Omega)} \quad \text{for every $v \in BV(\Omega)$},
		\end{equation}
		(We will always consider $\mathcal{S}_1$ to be the best constant for this embedding).
		\\Moreover, the embedding $BV(\Omega) \hookrightarrow L^p(\Omega)$ is compact for every $1 \le p < \frac{N}{N-1}$.
	\end{theorem}
	\medskip
	At this point, we review the chain rule for functions of bounded variation (\cite[Theorem 3.99]{AFP}).
	\begin{lemma}\label{lem c r}
		Let $v \in BV(\Omega)$ and let $\Phi: \mathbb{R} \to \mathbb{R}$ be a Lipschitz function. Then $w=\Phi(v) \in BV(\Omega)$ and
		\begin{equation*}
			Dw= \Phi'(\tilde{v})\tilde{D}v + \left(\Phi(v^e)-\Phi(v^i)\right)\nu_v \mathcal{H}^{N-1} \res J_v,
		\end{equation*}
		where $\tilde{D}v=D^av+D^cv$.
	\end{lemma}
	Observe that then
	\begin{equation}\label{ch wj}
		\tilde{D}w=\Phi'(\tilde{v})\tilde{D}v\quad \forall v \in DBV(\Omega).
	\end{equation}
	\medskip
Following \cite{GMP, BOPS}, to achieve our objectives, we propose the following function space.
	\begin{multline*}
		TBV(\Omega):=\left\{v: \Omega \to \mathbb{R} \quad \text{Lebesgue measurable}: \right.  \\ \left.  F(v^+), F(v^-)\in BV(\Omega) \quad\text{ for all } F \in W^{1,\infty}\left([0,\infty); [a,\infty)\right),\text{ where } a>0 \right\},
	\end{multline*}
	where $v^+:=\max\{0, v\}$ and $v^-:=\max\{0, -v\}$.
	It is possible to equivalently define $TBV(\Omega)$ as follows
	\begin{multline*}
		TBV(\Omega):=\left\{v: \Omega \to \mathbb{R} \quad \text{Lebesgue measurable}: \right. \\ \left. T_{a,b}(v^+), T_{a,b}(v^-)\in BV(\Omega), \,\,\, \text{for all $0<a<b< \infty$}\right\}.
	\end{multline*}
	A related space to $TBV(\Omega)$ is $GBV(\Omega)$, that can be found in \cite[Section 4.5]{AFP}. It is a space that shares many features with $BV(\Omega)$, as the existence of trace. We will follow the development that \cite{AFP} makes for $GBV(\Omega)$, adapting it to $TBV(\Omega)$.
	\\We introduce the concept of a trace on the boundary of $\Omega$ for functions belonging to $TBV(\Omega)$, this result is proven in \cite[Lemma 5.1]{GMP} for nonnegative functions.
	\begin{lemma}\label{lem trac TBV}
		Let us consider $0 \le v \in TBV(\Omega)$. Then, there exists a nonnegative function $v^\Omega \in L^1(\partial \Omega)$ such that
		$$\lim_{\rho \to 0} \frac{1}{|\Omega \cap B_\rho(x_0)|} \int_{\Omega \cap B_\rho(x_0)} |v(y)-v^\Omega(x_0)|\, dy=0 \quad \text{for $\mathcal{H}^{N-1}$-a.e. $x_0 \in \partial \Omega$.}$$
		Moreover,
		$$v^\Omega(x)= \lim_{\substack{a \to 0 \\ b \to \infty}}   \left(T_{a,b}(v)\right)^\Omega(x) \quad \text{for $\mathcal{H}^{N-1}$-a.e. in $x \in\partial \Omega$},$$
		and
		$$F(v^\Omega)=\left(F(v)\right)^\Omega \quad \text{for all $F \in W^{1,\infty}\left([0,\infty); [a,b]\right)$, for all $0<a<b<\infty$.}$$
	\end{lemma}
	
	\medskip
	As in the case of $GBV(\Omega)$, we cannot work with the concepts of approximate continuity or jump points.
	We now present the definition of the set of weak approximate jump points (for more details see \cite[Definition 4.28]{AFP}). Let $v \in L^1(\Omega)$, we define the upper and lower approximate limits of $v$ at $x$ as the following quantities, respectively
	$$	v^{\lor}(x):=\inf\{t\in \overline{\mathbb{R}}:\lim_{\rho \to 0}\rho^{-N}|\{v>t\}\cap B_{\rho}(x)|=0\},
	$$
	$$
	v^{\land}(x):=  \sup \bk\{t\in\overline{\mathbb{R}}:\lim_{\rho \to 0}\rho^{-N}|\{v<t\}\cap B_{\rho}(x)|=0\},$$ where $\overline{\mathbb{R}}:=\mathbb{R}\cup \left\{-\infty, \infty\right\}$.
	\\We emphasize that $v^\lor,v^\land: \Omega \to \overline{\mathbb{R}}$ are Borel measurable functions.
	We denote by $S_v^*$ the following set:
	$$
	S_v^* := \left\{ x \in \Omega \,:\, v^\land(x) < v^\lor(x) \right\}.
	$$
	Recall that
	$$
	v^\lor(x) = v^\land(x) = \tilde{v}(x) \quad \text{for all  $x \in \Omega \setminus S_v$},
	$$ which implies $S_v^* \subseteq S_v$.
	Furthermore, we define
	$$
	DTBV(\Omega) := \left\{ v \in TBV(\Omega) \,:\, \mathcal{H}^{N-1} \left( S_v^* \right) = 0 \right\}.$$
	We now wish to define the weak approximate jump set $J_v^*$ of a function $v \in L^1(\Omega)$; $x\in \Omega$ is a weak approximate jump point if there exists a unit vector $\nu_v^*(x)$ such that the weak approximate limit of the restriction of $v$ to the hyperplane $\left\{y \in \Omega \,:\, (y-x) \cdot \nu_v^*(x)>0\right\}$ is $v^\lor(x)$ and the weak approximate limit of the restriction of $v$ to the hyperplane $\left\{y \in \Omega \,:\, (y-x) \cdot \nu_v^*(x)<0\right\}$ is $v^\land(x)$.
	\\If $ v \in L^1(\Omega) $, it can be proven that (see \cite[Definition 4.30]{AFP} and subsequent comments)
	\begin{equation*}
		\begin{aligned}
			&J_v \subseteq J_v^*, \quad v^\lor(x)=\max \left\{v^e(x), v^i(x)\right\}, \quad v^\land(x)=\min \left\{v^e(x), v^i(x)\right\}, \\ &\nu_v^*(x)=\pm \nu_v(x), \quad \text{for all $x \in J_v$}.
		\end{aligned}
	\end{equation*}
	At this point, we can state a key result for functions in $ TBV(\Omega) $ that provides a characterization of the previously introduced sets (refer to \cite[Theorem 4.34]{AFP}).
	\begin{lemma}\label{salti in TBV}
		Let us consider a nonnegative function $v \in TBV(\Omega)$. Then
		\begin{enumerate}
			\item $S^*_v=\bigcup_{\substack{a \to 0 \\ b \to \infty}}  S_{T_{a,b}(v)}$ and $$v^\lor(x)=\lim_{\substack{a \to 0 \\ b \to \infty}}  (T_{a,b}(v))^\lor (x), \qquad v^\land(x)=\lim_{\substack{a \to 0 \\ b \to \infty}}  (T_{a,b}(v))^\land (x). $$
			\item $S^*_v$ is countably $\mathcal{H}^{N-1}$-rectifiable and $\mathcal{H}^{N-1}(S^*_v \setminus J^*_v)=0$.
		\end{enumerate}
	\end{lemma}
	We conclude this overview of $TBV$-functions by explicitly noting that, as a consequence of the coarea formula (see \cite[Theorem 3.40]{AFP}), the sets $\{v > a\}$ and $\{v < -a\}$ possess finite perimeter for a.e. $a > 0$, provided that $v \in TBV(\Omega)$. As a result, the characteristic functions $\chi_{\{a < v < b\}}$ and $\chi_{\{-b < v < -a\}}$ belong to $BV(\Omega)$ for a.e. choice of $0<a<b<\infty$.
	\section{Generalized pairings having a divergence--measure field}\label{sec3}
	In this paper, we need to handle ``products'' of non-bounded divergence--measure vector fields and gradients of certain functions of bounded variation. This section is devoted to give sense to these pairings and to show some useful features.
	
	\subsection{The Anzellotti-Chen-Frid pairing}
	We begin this section by recalling the theory of pairing introduced by Anzellotti in \cite{A} and later developed by Chen and Frid in \cite{CF}.
	\\We define the following space
	$$\DM(\Omega):=\left\{z \in L^\infty(\Omega)^N \,:\, \operatorname{div}z \in \mathcal{M}(\Omega)\right\}.$$
	For any vector field $ z \in \DM(\Omega) $, there exists a generalized notion of trace on $ \partial \Omega $ of its normal component. Indeed, the theorem asserts that there exists a linear operator $ [\cdot, \nu]: \DM(\Omega) \to L^\infty(\partial \Omega) $ such that
	\begin{equation*} 
		\|[z,\nu]\|_{L^\infty(\partial\Omega)} \le \|z\|_{L^\infty(\Omega)^N},
	\end{equation*}
	and if $z\ \in C^1(\overline{\Omega})^N$, it can prove that
	$$[z,\nu]=z(x)\cdot \nu(x) \quad \text{for $\mathcal{H}^{N-1}$-a.e. $x \in \partial\Omega$.}$$
	Furthermore, \cite[Proposition 3.1]{CF} asserts that $ \operatorname{div}z $ is absolutely continuous with respect to $ \mathcal{H}^{N-1} $.
	\\At this point, we introduce the following distribution $\left(z, Dv\right): C^\infty_c(\Omega) \to \mathbb{R}$
	\begin{equation}\label{def pair}
		\langle \left(z, Dv\right), \varphi\rangle := - \into v^* \varphi \operatorname{div}z - \into vz \cdot \nabla \varphi,
	\end{equation}
	which is well-defined if we take $ z \in \DM(\Omega) $ and $ v \in BV(\Omega) \cap L^\infty(\Omega) $, as stated in \cite[Appendix A]{MST} or  \cite[Section 5]{C}.
	Then the pairing $ (z, Dv) $ is a Radon measure with finite total variation, it satisfies the inequality
	\begin{equation} \label{pair bel mis}
		|\langle (z, Dv), \varphi \rangle| \leq \|\varphi\|_{L^\infty(A)} \|z\|_{L^\infty(A)^N} \int_A |Dv|\,,
	\end{equation}
	for every open set $ A\Subset \Omega $ and every $ \varphi \in C_c^1(A) $. In particular,  $ |(z, Dv)|$ is absolutely continuous with respect to $|Dv| $.
	Therefore, by using the Radon-Nikodým Theorem, we will obtain:
	\begin{equation}\label{pair con der teta}
		\left(z, Dv\right)= \theta\left(z, Dv,x\right) |Dv| \quad \text{as measures in $\Omega$},
	\end{equation}
	where $\theta(z, Dv, x)$ is the Radon-Nikodým derivative of $(z,Dv)$ with respect to $|Dv|$.
	As for the Radon-Nikodým derivative, we can establish the following result (we invoke \cite[Proposition 4.5 (iii)]{CDC}; see also \cite[Proposition 2.8]{A} and \cite[Proposition 2.7]{LS}).
	\begin{lemma}\label{prop 4.5 cdc}
		Let us choose $z \in \DM(\Omega)$, $v \in BV(\Omega) \cap L^\infty(\Omega)$ and $\Phi: \mathbb{R} \to \mathbb{R}$ a nondecreasing Lipschitz function. Then
		\begin{equation}\label{eq 4.5 cdc}
			\theta\left(z, D\Phi(v), x\right)=\theta\left(z, Dv, x\right) \quad \text{for $|D\Phi(v)|$-a.e. $x \in \Omega$.}
		\end{equation}
	\end{lemma}
	For our purposes, it is useful to state the following result (see \cite[Lemma 5.4, Lemma 5.6]{C}).
	\begin{lemma}\label{lem uz in DM}
		Let $z \in \DM(\Omega)$ and let $v \in BV(\Omega) \cap L^\infty(\Omega)$, then it holds \begin{equation}\label{leib rule dm infty}
			\left(z,Dv\right)=-v^* \operatorname{div}z+\operatorname{div}\left(\displaystyle vz\right) \quad \text{as measures in $\Omega$.}
		\end{equation} Moreover,  $vz \in \DM(\Omega)$ and it holds \begin{equation}\label{v fuori v dentro}
			[v z,\nu]=v [z,\nu] \quad \text{for $\mathcal{H}^{N-1}$-a.e. $x \in \partial \Omega$.}
		\end{equation}.
	\end{lemma}
	Now, let us present the following Lemma proven in \cite[Proposition 2.3]{MS} (see also \cite[Lemma 2.6]{GMP}).
	\begin{lemma}
		Let us consider $z \in \DM(\Omega)$ and suppose $u, v \in DBV(\Omega) \cap L^\infty(\Omega)$. Then
		\begin{equation}\label{uscire sx}
			\left(uz, Dv\right)=u\left(z, Dv\right) \quad \text{as measures in $\Omega$.}
		\end{equation}
	\end{lemma}
	At this point, we provide the statement of the Gauss-Green formula in this more generalized context, as proven in \cite{CF} (see also \cite[Theorem A.1]{MST} and \cite[Theorem 5.3]{C}).
	\begin{lemma}
		Let $z \in \DM(\Omega)$ and $v \in BV(\Omega) \cap L^\infty(\Omega)$. Then, the following generalized Gauss-Green formula holds
		\begin{equation}\label{GG form}
			\int_{\Omega} v^* \, \operatorname{div} z + \int_{\Omega}(z,Dv) = \int_{\partial \Omega} v [z,\nu] \, \ensuremath{\mathrm d}\mathcal{H}^{N-1}.
		\end{equation}
	\end{lemma}
	
	\medskip
	It is also important to note that the normal trace $ [z, \xi]^\pm $ of a vector field $ z \in \DM(\Omega) $ on an oriented $ C^1 $-hypersurface $ \xi \subset \Omega $ can be defined as
	$$
	[z, \nu_\xi]^\pm := [z, \nu_{\omega^\pm}], $$
	where $ \omega^\pm \Subset \Omega $ are open $ C^1 $-domains such that $ \xi \subset \partial \omega^\pm $ and $ \nu_{\omega^\pm} = \pm \nu_\xi $. This definition does not depend on the specific choice of $ \omega^\pm $, except on a set $ \mathcal{H}^{N-1}$-negligible. Moreover, as established in \cite[Proposition 3.4]{ACM}, the following relation holds
	\begin{equation}\label{t div su ip}
		\left(\operatorname{div} z\right)\res \xi = \left(\left[z, \nu_\xi\right]^+ - \left[z, \nu_\xi\right]^-\right) \mathcal{H}^{N-1} \res \xi.
	\end{equation}
	This concept can be generalized by localization to oriented countably $ \mathcal{H}^{N-1} $-rectifiable sets, enabling an extension of the above formula  (see \cite[Lemma 2.4]{GMP}).
	
	As a consequence of Lemma \ref{lem uz in DM}, the following result can also be derived (see \cite[Lemma 2.5]{GMP}):
	
	\begin{lemma} 
		Let $ v \in BV(\Omega) \cap L^\infty(\Omega) $ and $ z \in \mathcal{DM}^\infty(\Omega) $. Then, for every $\xi \subset \Omega$ regular hypersurface, we have
		\begin{equation}\label{saltando sul bordo}
			\left[v z, \nu_\xi\right]^\pm = v^\pm \left[z, \nu_\xi\right]^\pm \quad \text{for $ \mathcal{H}^{N-1} $-a.e. on $ \xi $}.
		\end{equation}
	\end{lemma}
	\subsection{Pairing with unbounded vector field} 
	To achieve our goal, it is necessary to present the concept of pairing in the setting of an unbounded vector field \( z \). To this end, we follow the generalization of this tool introduced in \cite{CDCS}; in our specific case, we will adapt the notions already seen to the following vector fields $$\mathcal{DM}^p(\Omega):=\left\{z \in L^p(\Omega)^N \,:\, \operatorname{div}z \in \mathcal{M}(\Omega)\right\} \quad \text{with $1 \le p < \infty$} .$$
	\\Let us begin by providing the concept of $\lambda$-representative of a Borel function. Given $v \in \mathcal{B}(\Omega)$ and $\lambda \in \mathcal{B}(\Omega; [0,1])$, we define the $\lambda$-representative $v^\lambda: \Omega \to \overline{\mathbb{R}} $ as follows
	\begin{equation*}
		v^\lambda(x):=
		\begin{cases}
			(1-\lambda(x))v^\land(x)+\lambda(x)v^\lor(x) & \text{if $x \in \Omega \setminus Z_v$,}\\
			\infty & \text{if $x \in Z_v$ and $\lambda(x)>\frac{1}{2}$,}\\
			0 & \text{if $x\in Z_v$ and $\lambda(x)=\frac{1}{2}$,}\\
			-\infty & \text{if $x \in Z_v$ and $\lambda(x)<\frac{1}{2}$,}
		\end{cases}
	\end{equation*}
	where $Z_v:=\left\{x \in \Omega \,:\, v^\lor(x)=\infty \quad \text{and} \quad v^\land(x)=-\infty\right\}$ and $v^\land, v^\lor$ are the function introduced in \ref{subs BV}. In particular, when $\lambda \equiv 1$ and $\lambda \equiv 0$, we have $v^1:=v^\lor$ and $v^0:=v^\land$. Moreover, if $v \in L^1(\Omega)$ then $Z_v \subseteq S_v \setminus J_v$ and if $v \in L^\infty(\Omega)$, then $Z_v=\emptyset$, which implies that, whatever $\lambda$ is, $$v^\lambda(x)=v(x) \quad \text{for $\mathcal{L}^N$-a.e. $x\in \Omega$},$$ because $|S_v^*|=0$.
	\\In the case $v \in BV(\Omega)$, we obtain $\mathcal{H}^{N-1}(Z_v) \le \mathcal{H}^{N-1}(S_v \setminus J_v)=0$, as stated in subsection \ref{subs BV}. As a consequence, if $\lambda(x)=\frac12$ on $J_v$, then  $v^{\lambda}(x)=v^*(x)$ for every $x \in \Omega\setminus (S_v \setminus J_v)$ (for more details, we refer to \cite[p.11]{CDCS}).
	\\Given $z \in \mathcal{DM}^1(\Omega)$ and $\lambda \in \mathcal{B}(\Omega; [0,1])$, let us consider the following set  $$\mathcal{X}^{z, \lambda}(\Omega):=\left\{v \in \mathcal{B}(\Omega) : v^\lambda \in L^1(\Omega, |z|)\cap L^1(\Omega, |\operatorname{div}z|)\right\},$$ we underline that $\mathcal{X}^{z, \lambda}(\Omega)$ is not a vector space since the $\lambda$-representative of the sum is not the sum of $\lambda$-representatives (see \cite[Remark 7.5]{CDCS}).
	\\We provide the definition of pairing in this generalized context, as given in \cite[Definition 3.1]{CDCS}.
	\begin{defin}
		Let $z \in \mathcal{DM}^1(\Omega)$. If $\lambda \in \mathcal{B}(\Omega; [0,1])$ and $v \in \mathcal{X}^{z,\lambda}(\Omega)$, then the distribution $(z,Dv)_\lambda:C^\infty_c(\Omega) \to \mathbb{R}$ \begin{equation*}
			\langle \left(z, Dv\right)_\lambda, \varphi\rangle := - \into v^\lambda \varphi \operatorname{div}z - \into vz \cdot \nabla \varphi ,
		\end{equation*}
		is well defined. (We stress that $v^\lambda=v$ for $\mathcal L^N$-almost all $x\in \Omega$, so that we may write just $v$ in the second integral).
	\end{defin}
	It is important to emphasize that for a function $v \in \mathcal{X}^{z,\lambda}(\Omega)$, the derivative $Dv$ is not, in general, a Radon measure, but it is always well-defined as a distribution of order $1$.
	\begin{remark}
		In the particular case where $v \in W^{1,1}(\Omega) \cap L^\infty(\Omega)$, we will have (for more details see \cite[Proposition 3.15 (i)]{CDCS}) $$\left(z,Dv\right)_\lambda= z \cdot \nabla v \, \mathcal{L}^N.$$
	\end{remark}
	Drawing an analogy with the classical theory of functions of bounded variation, we present a version of $BV$-type classes defined in terms of the variation induced by this generalized $\lambda$-pairing (see \cite[Definition 3.3]{CDCS}).
	\begin{defin}
		Let us give $z \in \mathcal{DM}^1(\Omega)$ and $\lambda \in \mathcal{B}(\Omega; [0,1])$, we define the set
		$$BV^{z,\lambda}(\Omega):=\left\{v \in \mathcal{X}^{z,\lambda}(\Omega):(z,Dv)_\lambda \in \mathcal{M}(\Omega)\right\}.$$
	\end{defin}
	Also, $BV^{z,\lambda}(\Omega)$ is not a vector space; for further details, we refer to \cite[Remark 3.4]{CDCS}.
	\\To extend Lemma \ref{lem uz in DM} to this new context, we present the following result (see \cite[Proposition 3.5]{CDCS}).
	\begin{lemma}\label{lemma per mis l1}
		Let $z \in \mathcal{DM}^1(\Omega)$ and let $\lambda \in \mathcal{B}(\Omega; [0,1])$. Then   $v\in BV^{z,\lambda}(\Omega)$ if and only if $v \in \mathcal{X}^{z,\lambda}(\Omega)$ and $\operatorname{div}\left(\displaystyle vz\right) \in \mathcal{M}(\Omega)$. Moreover, we get		
		\begin{equation}\label{ug pair con div lambda}
			\left(z,Dv\right)_\lambda=-v^\lambda \operatorname{div}z+\operatorname{div}\left(\displaystyle vz\right) \quad \text{as measures in $\Omega$.}
		\end{equation}
	\end{lemma}
	The next result is a special case of the previous ones. It is worth stating here, since we will use it in what follows (see \cite[Proposition 3.19]{CDCS}).
	\begin{proposition}\label{prop z e div in l1}
		Let us take $z \in \mathcal{DM}^1(\Omega)$ such that $|\operatorname{div}z| \ll \mathcal{L}^N$. Then $\mathcal{X}^{z,\lambda}(\Omega)$ does not depend on $\lambda$. In other words,  for every $\lambda \in \mathcal{B}(\Omega; [0,1])$, it holds
		$$\mathcal{X}^{z,\lambda}(\Omega)=\mathcal{X}^{z}(\Omega)=\left\{v \in \mathcal{B}(\Omega)\,:\, v\in L^1(\Omega, |z|)\cap L^1(\Omega, |\operatorname{div}z|)\right\},$$
		and if $v \in \mathcal{X}^z(\Omega)$, then the distribution $(z,Dv):C^\infty_c(\Omega) \to \mathbb{R}$ given by
		\begin{equation}\label{pair l1}
			\langle \left(z, Dv\right), \varphi\rangle := - \into v \varphi \operatorname{div}z - \into vz \cdot \nabla \varphi ,
		\end{equation}
		is well defined.
		Moreover, it follows
		$$BV^{z,\lambda}(\Omega)=BV^z(\Omega)=\left\{v \in \mathcal{X}^z(\Omega)\,:\,  \left(z, Dv \right)\in \mathcal{M}(\Omega)\right\},$$ and, for every $v\in BV^z(\Omega)$, we have
		\begin{equation}\label{ug pair con div l1}
			\left(z, Dv\right)_\lambda=(z,Dv)=-v\operatorname{div}z+\operatorname{div}\left(\displaystyle vz\right) \quad \text{as measures in $\Omega$.}
		\end{equation}
	\end{proposition}
	
	We recall that when $z \in \mathcal{DM}^\infty(\Omega)$, $\lambda \in \mathcal{B}(\Omega; [0,1])$, and $v \in BV^{z,\lambda}(\Omega)\cap L^\infty(\Omega)$, the generalized pairing $(z, Dv)_{\lambda}$ satisfies the following property:
	$$|(z,Dv)_{\lambda}| \ll \mathcal{H}^{N-1}, $$
	as stated in \cite[Proposition 3.14(ii)]{CDCS}.

	\begin{remark}\label{abs-cont-Haus}
	We point out that this property cannot be extended to a more general context. Indeed, in \cite[Proposition 3.15(ii)]{CDCS}, the authors demonstrate that if $z \in \mathcal{DM}^p(\Omega)$ with $p \ge \frac{N}{N-1}$ and $v \in BV^{z,\lambda}(\Omega)\cap L^\infty(\Omega)$, then
	$$|(z,Dv)_\lambda|(B)=0 \quad  \text{for all sets $B\subset \Omega$ which are $\sigma$-finite with respect to $\mathcal{H}^{N-\frac{p}{p-1}}$} .$$
	Moreover, if $ 1 \le p < \frac{N}{N-1}$, then the $\lambda$-pairing is not, in general, absolutely continuous with respect to any Hausdorff measure, because the singularities of $\operatorname{div}z$ can be arbitrary (for further details, see \cite[Remark 3.17]{CDCS}).
	
	This property has an important consequence concerning the notion of the trace of a vector field $z \in \mathcal{DM}^p(\Omega)$, with $1 \le p < \frac{N}{N-1}$, which a priori can only be a distribution. As can be seen in \cite[p.15]{CDCS}, given $\omega \Subset \Omega$ with a regular boundary (for such sets, $\chi_\omega \in BV^{z,\lambda}(\Omega)$, for every $\lambda \in \mathcal{B}(\Omega; [0,1])$; for more details, we refer to \cite[Section 7]{CDCS}), we can generalize the concepts of the interior and exterior traces, respectively, as follows:
	$$(z, D\chi_\omega)_0 \quad \text{and} \quad (z, D\chi_\omega)_1.$$
	
	But since $(z, D \chi_\omega)_0$ and $(z, D \chi_\omega)_1$ are not, in general, absolutely continuous with respect to any Hausdorff measure, and in particular not with respect to $\mathcal{H}^{N-1}$, we cannot repeat the same arguments as in \cite[Proposition 4.7]{CCDCM}, where the authors prove that if $z \in \mathcal{DM}^\infty(\Omega)$, then $[z,\nu_\omega]^+$ and $[z,\nu_\omega]^-$ are, respectively, the densities of the measures $(z, D\chi_\omega)_0$ and $(z, D\chi_\omega)_1$ with respect to the measure $|D\chi_\omega|$.
	\end{remark}
	\bk
	
	\section{Statement of main result}\label{sec4}
	Our aim is to study the following Dirichlet problem for the relativistic transparent media operator
	\begin{equation}\label{Prob tm1}
		\begin{cases}
			-\operatorname{div}\left(\displaystyle u^m \frac{Du}{|Du|}\right)=f & \text{in $\Omega$,}\\
			u=0 & \text{on $\partial \Omega$,}
		\end{cases}
	\end{equation}
	where $m>0$ and $0\le f \in L^1(\Omega)$.
	Let us introduce the notion of solution to problem \eqref{Prob tm1}.
	\begin{defin}\label{def sol}
		Let $m>0$ and $0\le f \in L^1(\Omega)$. A nonnegative $u \in DTBV(\Omega)$ satisfying 	$u^m\in L^{\frac {N}{N-1}, \infty}(\Omega)$, $T_b(u)^{m+1} \in BV(\Omega)$ for a.e. $b>0$ is solution  to  problem \eqref{Prob tm1} if there exists a vector field $w \in L^\infty(\Omega)^N$ with $\|w\|_{L^\infty(\Omega)^N}\le 1$ such that $z:=u^m w$ belongs to $\mathcal{DM}^1(\Omega)$ and it has a weak trace $[z,\nu] \in L^1(\partial\Omega)$, and the following conditions hold
		\begin{equation}\label{sol 1}
			-\operatorname{div}z=f \quad \text{as distributions in $\Omega$,}
		\end{equation}
		\begin{equation}\label{sol 2}
			\left(z, DT_{a,b}(u)\right)= \frac{1}{m+1}\left|DT_{a,b}(u)^{m+1}\right| \quad \text{as measures in $\Omega$,}
		\end{equation}
		for a.e. $0<a<b<\infty$, and
		\begin{equation}\label{sol 3}
			[z,\nu]=-(u^\Omega)^m \quad \text{$\mathcal{H}^{N-1}$ on $\partial \Omega \cap \{u>0\}$.}
		\end{equation}
	\end{defin}
	\begin{remark}
		Let us take a moment to offer some clarifications regarding Definition \ref{def sol}. Equation \eqref{sol 1} reveals how the vector field $w$ serves as a weak interpretation of the quotient $\frac{Du}{|Du|}$, while $z$ reflects the expression  $u^m \frac{Du}{|Du|}$. This relationship is further explained by the subsequent identity \eqref{sol 2}.
		\\Moreover, we emphasize the significance of specifying that the vector field $z \in \mathcal{DM}^1(\Omega)$ possesses a weak trace $[z, \nu] \in L^1(\partial \Omega)$. Indeed, in the event that the vector field is unbounded, one can only guarantee that its trace exists as a well-defined distribution (for further details, see Remark \ref{abs-cont-Haus}). It is also possible to define a notion of trace for functions $u \in DTBV(\Omega)$, as shown in Lemma \ref{lem trac TBV}.
		\\Thanks to these considerations, the boundary condition \eqref{sol 3} is well-posed for $\mathcal{H}^{N-1}$-a.e. $x \in \partial \Omega$. As a consequence, the trace $(u^\Omega)^m\in L^1(\partial\Omega)$.
	\end{remark}
	Now we state the main result of this paper.
	\begin{theorem}\label{teo tras 1}
		Let $m>0$ and $0 \le f \in L^1(\Omega)$. Then, there exists a nonnegative function $u \in DTBV(\Omega)$ which is a solution to the problem \eqref{Prob tm1} in the sense of Definition \ref{def sol}.
	\end{theorem}
	\begin{remark}
	We highlight the key fact that, in the case $m > 0$, the solution obtained is finite $\mathcal{L}^N$-a.e., as established in Remark \ref{rem u fqo}. This stands in sharp contrast to what occurs in the case of the 1-Laplacian operator (i.e., the case $m=0$), where such a finite value is not generally guaranteed (see, for instance, \cite{MST}).
	\end{remark}
	
	\section{Proof of main result}\label{sec5}
	
	\subsection{Approximating problems and main estimates}\label{subs aprx p}
	We present the following approximation scheme to address our problem
	\begin{equation}\label{Prob trans tronc}
		\begin{cases}
			-\operatorname{div}\left(\displaystyle u_n^m \frac{Du_n}{|Du_n|}\right)=f_n & \text{in $\Omega$,}\\
			u_n=0 & \text{on $\partial \Omega$,}
		\end{cases}
	\end{equation}
	with $f_n=T_n(f)$ ($T_n(s)$ is the truncation function defined in \eqref{T_k}).
	\medskip
	\\In \cite[Theorem 3.3]{BOPS}, it is proven that, for every fixed $n \in \mathbb{N}$, there exists a nonnegative solution $u_n \in DTBV(\Omega) \cap L^\infty(\Omega)$ such that $u_n^{m+1} \in BV(\Omega)$, and an associated vector field $w_n \in L^\infty(\Omega)^N$ with $\|w_n\|_{L^\infty(\Omega)^N} \leq 1$ and $z_n :=u_n^m w_n \in \DM(\Omega)$, which satisfy the following conditions
	\begin{equation}\label{sol 1 n}
		-\operatorname{div}z_n=f_n \quad \text{as distributions in $\Omega$,}
	\end{equation}
	\begin{equation}\label{sol 2 n}
		\left(z_n, DT_{a,b}(u_n)\right)=\frac{1}{m+1}\left|DT_{a,b}(u_n)^{m+1}\right| \quad \text{as measures in $\Omega$, }
	\end{equation}
	for a.e. $0<a<b \le \infty$, and
	\begin{equation}\label{sol 3 n}
		[z_n,\nu]=-(u_n^\Omega)^m \quad \text{$\mathcal{H}^{N-1}$-a.e. on $\partial\Omega \cap \{u_n>0\}$.}
	\end{equation}
	\medskip
	\begin{remark} 
		We point out that the identity
		\[\left(z_n, DT_{a,b}(u_n)\right)=\left(z_n \chi_{\{a<u_n<b\}}, DT_{a,b}(u_n)\right)\quad \text{as measures in $\Omega$,}\]
		holds for almost all $0<a<b \le \infty$.
		Furthermore, the preceding equality allows us to deduce an alternative version of \eqref{sol 2 n}
		\begin{equation*}
			\left(w_n \chi_{\{a<u_n<b\}},DT_{a,b}(u_n) \right)= |DT_{a,b}(u_n)| \quad \text{as measures in $\Omega$,}
		\end{equation*}
	for a.e. $0<a<b \le \infty$, as a consequence of \eqref{uscire sx} and the chain rule property \eqref{ch wj}.
	\end{remark}
	\begin{remark}\label{oss-2}
	The relationship $[z_n,\nu]=u_n^m[w_n,\nu]$ $\mathcal{H}^{N-1}$-a.e on $\partial \Omega$  would typically be expected. However, in the absence of regularity conditions (specifically, that $w$ is a divergence-measure field or that $u\in BV(\Omega)$), we cannot assert its validity.
		\end{remark}
	\begin{remark} 
		Having in mind Green's formula \eqref{GG form}, it holds
		\begin{equation}\label{exp f test n}
			\int_\Omega (z_n,D\psi)-\int_{\partial \Omega} \psi [z_n, \nu] \, \mathrm{d}\mathcal{H}^{N-1} =\int_\Omega f_n \psi \quad \text{for all $\psi \in BV(\Omega) \cap L^\infty(\Omega)$,}
		\end{equation}
		which is the weak formulation of problem \eqref{Prob trans tronc}.
	\end{remark}
	\medskip
	Initially, we present some estimates on $u_n$ that ensure the existence of a limit function $u$, which serves as our candidate to solution.
	\begin{lemma}
		Assume $m>0$ and $0\le f \in L^1(\Omega)$. Let us consider $u_n$ solution to the problem \eqref{Prob trans tronc}. Then for a.e. $b>0$, they hold
		\begin{equation}\label{stime T_b}
			\|T_b(u_n)^{m+1}\|_{BV(\Omega)}\le (m+1)\|f\|_{L^1(\Omega)} \, b,
		\end{equation}
		and
		\begin{equation}\label{stime Marc}
			[u_n^m]_{L^{{\frac{N}{N-1}}, \infty}(\Omega)}\le \left( (m+1)\mathcal{S}_1 \|f\|_{L^1(\Omega)} \right)^{\frac{N}{N-1}}.
		\end{equation}
	\end{lemma}
	\begin{proof}
		We begin by proving \eqref{stime T_b}.
		We choose $T_{a,b}(u_n)$ with $0<a<b<\infty$ as a test function in \eqref{exp f test n}, we gain
		\begin{equation}\label{prob ap con tab}
			\int_\Omega \left(z_n, DT_{a,b}(u_n)\right)- \int_{\partial\Omega} T_{a,b}(u_n) [z_n, \nu] \, \mathrm{d}\mathcal{H}^{N-1}=\int_\Omega f_n T_{a,b}(u_n).
		\end{equation}
	We remind the reader that \cite[Lemma 5.8]{BOPS} provided the relation $T_{a,b}(u_n)^m \le - [z_n,\nu] $. So, by this fact and \eqref{sol 2 n}, \eqref{prob ap con tab} becomes
		\begin{equation*}\label{tab con sol 2}
			\frac{1}{m+1}\int_\Omega \left|DT_{a,b}(u_n)^{m+1}\right|+ \frac{1}{m+1}\int_{\partial\Omega} T_{a,b}(u_n)^{m+1} \, \mathrm{d}\mathcal{H}^{N-1} \le \int_\Omega f_n T_{a,b}(u_n),
		\end{equation*}
		which implies that
		\begin{equation*}
			\|T_{a,b}(u_n)^{m+1}\|_{BV(\Omega)}\le (m+1) \|f\|_{L^1(\Omega)} \, b.
		\end{equation*}
		We observe that the right hand side is independent of $a > 0$. Thanks to the fact that $u_n^{m+1} \in BV(\Omega)$, taking the limit as $a$ tends to $0$, we derive \eqref{stime T_b}.
		\medskip
		\\We now turn to demonstrate the estimates \eqref{stime Marc}. Applying Sobolev's embedding \eqref{sob emb} to \eqref{stime T_b}, it yields
		\begin{equation*}
			|\{u_n>b\}|b^{(m+1)\frac{N}{N-1}}\le\into |T_b(u_n)^{m+1}|^{\frac{N}{N-1}} \le \left(\mathcal{S}_1(m+1)  \|f\|_{L^1(\Omega)} \, b\right)^{\frac{N}{N-1}}.
		\end{equation*}
		Setting $h^{\frac{1}{m}} = b$ and noting that $\{u_n^m > h\} = \{u_n > h^{\frac{1}{m}}\}$, we obtain
		\begin{equation*}
			|\{u_n^m>h\}|\le \frac{\left(\mathcal{S}_1(m+1)\|f\|_{L^1(\Omega)}\right)^{\frac{N}{N-1}}}{h^{\frac{N}{N-1}}},
		\end{equation*}
		which implies \eqref{stime Marc} and this concludes the proof.
	\end{proof}
	\medskip
	In the next result we prove the existence of a nonnegative function $u$ which is our candidate to solution.
	\begin{corollary}\label{cor ext u}
		Assume $m>0$ and $0 \le f \in L^1(\Omega)$. Let us consider $u_n$ solution to the problem \eqref{Prob trans tronc}. Then there exists a nonnegative measurable function $u$ such that $u_n^m \to u^m$ strongly in $L^q(\Omega)$ for every $1\le q < \frac{N}{N-1}$ and $T_{a,b}(u_n) \to T_{a,b}(u)$ strongly in $L^q(\Omega)$ for every $1 \le q < \infty$ and $DT_{a,b}(u_n) \rightharpoonup DT_{a,b}(u)$ $*$-weakly as measures, for a.e. $0<a<b<\infty$.
	\end{corollary}
	\begin{proof}
		Initially, we show the existence of a nonnegative measurable function $u$.
		By \eqref{stime T_b} and the compact properties of $BV(\Omega)$ \eqref{sob emb}, for a.e. $b>0$ there exists $g_b^{m+1} \in BV(\Omega)$ such that
		\begin{equation}\label{conv delle T_b}
			T_b(u_n)^{m+1} \to g_b^{m+1} \quad \text{ strongly in $L^q(\Omega)$ for every $1 \le q < \frac{N}{N-1}$ and a.e.}
		\end{equation}
		\medskip Now we define a function $u: \Omega \to \mathbb{R}$ such that
		\begin{equation*}
			u(x):=
			\begin{cases}
				g_b(x) & \text{ if $x \in  \{g_b(x)<b\}$, for some  $b>0$,} \\
				\infty & \text{ if  $x\notin \cup_{b>0} \{g_b(x)<b\}\,$.}
			\end{cases}
		\end{equation*}
		We underline that the function $u$ is well-defined, as for each $x \in \Omega$, the value of $u(x)$ is independent of the choice of $b$. In fact, by choosing $0 < c < d < \infty$, we have $T_c(u_n(x)) = T_d(u_n(x))$ for every $x \in \{u_n < c\}$. Consequently,
		$$g_d(x)^{m+1} \stackrel{\eqref{conv delle T_b}}{=}\lim_{n \to \infty} T_d(u_n(x))^{m+1} = \lim_{n \to \infty} T_c(u_n(x))^{m+1} \stackrel{\eqref{conv delle T_b}}{=} g_c(x)^{m+1}.$$
		We stress that $u$ is a Lebesgue measurable function, since it is the pointwise limit of measurable functions.
		\medskip
		\\
		It follows from \eqref{conv delle T_b} and the definition of $u$  that
		\begin{equation}\label{conv TAB}
			T_{a,b}(u_n) \to T_{a,b}(u) \quad \text{strongly in $L^q(\Omega)$ for every $1\le q <\infty$,}
		\end{equation}
		for a.e. $0<a<b<\infty$. Recalling 	 Theorem \ref{embedding}, we can also assure that
		\begin{equation*}\label{conv DTAB}
			DT_{a,b}(u_n) \rightharpoonup DT_{a,b}(u) \quad \text{ $*$-weakly as measures in $\Omega$, for a.e. $0<a<b<\infty$.}
		\end{equation*}
		\medskip
		Moreover, as a consequence of \eqref{conv delle T_b}, we have
		\begin{equation}\label{qo di um}
			u_n^m \to u^m \quad \text{a.e. $x \in \Omega$.}
		\end{equation}
		It just remains to check  that
		\begin{equation}\label{conv forte u^m}
			u_n^m \to u^m \quad \text{strongly in $L^q(\Omega)$ for every $1\le q < \frac{N}{N-1}$.}
		\end{equation}
		First notice that, owing to \eqref{stime Marc}, the sequence $u_n^m$ is bounded in $L^q(\Omega)$ for every $1 \le q<\frac N{N-1}$. Thus, Hölder's inequality implies that $u_n^m$ is equi-integrable. Invoking Vitali's Theorem, we deduce that
		\begin{equation}\label{conv in L^1}
			u_n^m \to u^m \quad \text{strongly in $L^1(\Omega)$.}
		\end{equation}
		Finally, given $1 \le q<\frac N{N-1}$, choose $r\in (q,\frac{N}{N-1})$. Since  $u_n^m$ is bounded in $L^r(\Omega)$ and \eqref{conv in L^1} holds, it follows from the interpolation inequality that
		\begin{equation*}
			u_n^m \to u^m \quad \text{strongly in $L^q(\Omega)$,}
		\end{equation*}
		as desired.
	\end{proof}
	\begin{remark}\label{rem u fqo}
		From \eqref{stime Marc} and \eqref{qo di um}, we can state that
		$$|\{u > k\}| \le \frac{\left(\mathcal{S}_1(m+1)\|f\|_{L^1(\Omega)}\right)^{\frac{N}{N-1}}}{k^{m \frac{N}{N-1}}},$$
		which implies that $u \in L^{\frac{mN}{N-1}, \infty}(\Omega)$, and consequently,
		\begin{equation}\label{u fqo}
			\left|\left\{u = \infty\right\}\right| = 0.
		\end{equation}
	\end{remark}
	
	\subsection{Existence of the vector fields and condition (\ref{sol 1})}
	In the following results, we demonstrate the existence of a vector field $w \in L^\infty(\Omega)^N$ which will play the role of the quotient $\frac{Du}{|Du|}$.
	\begin{lemma}\label{lem ext z}
		Assume $m>0$ and $0 \le f \in L^1(\Omega)$. Let $u_n$ be a solution to the problem \eqref{Prob trans tronc} having the associated vector field $w_n$. Then, there exists a vector field $w \in L^\infty(\Omega)^N$ with $\|w\|_{L^\infty(\Omega)^N} \le 1$ such that $z:=u^m w$ satisfies \eqref{sol 1}.
		\\Moreover, $z \in \mathcal{DM}^1(\Omega)$ and $\operatorname{div}z$ is absolutely continuous with respect to $\mathcal{L}^N$.
	\end{lemma}
	\begin{proof}
		To show the existence of the vector field $w$, let us consider the sequence of vector field $w_n$. Since $ \|w_n\|_{L^\infty(\Omega)^N} \leq 1 $, there exists a vector field $ w \in L^\infty(\Omega)^N $ such that $ w_n \rightharpoonup w $ $*$-weakly in $L^\infty(\Omega)^N$ and so $ \|w\|_{L^\infty(\Omega)^N} \leq 1 $. Set $z:=u^mw \in L^{\frac{N}{N-1}, \infty}(\Omega)^N$.
		Thanks to \eqref{conv forte u^m}, we can assert that
		\begin{equation}\label{con deb zn to z}
			z_n \rightharpoonup z \quad \text{weakly in $L^q(\Omega)^N$, for every $ 1 \le q < \frac{N}{N-1}$.}
		\end{equation}
		\medskip
		Now, we prove \eqref{sol 1}, choosing $0\le \varphi \in C^1_c(\Omega)$ as a test function in \eqref{sol 1 n}, we get
		\begin{equation*}
			\into z_n \cdot \nabla \varphi = \into f_n \varphi.
		\end{equation*}
		Letting $n$ go to infinity, the left hand side passes to limit by \eqref{con deb zn to z}, while the right hand side passes due to Lebesgue's Theorem, so we get \eqref{sol 1} and this concludes the proof.
	\end{proof}
	
	\subsection{A first inequality in condition (\ref{sol 2})}\label{sub 5.3}
	To establish equation \eqref{sol 2}, we begin by proving an inequality that serves as its preliminary form.
	\begin{lemma}
		Assume $m>0$ and $0 \le f \in L^1(\Omega)$. Let $u_n$ be a solution to the problem \eqref{Prob trans tronc} with associated vector field $w_n$. Then, for a.e. $0<a<b<\infty$, the distribution  $\left(z, DT_{a,b}(u)\right)$ is a Radon measure which satisfies
		\begin{equation}\label{dis per sol 2}
			\frac{1}{m+1} \left|DT_{a,b}(u)^{m+1}\right| \le \left(z, DT_{a,b}(u)\right) \quad \text{as measures in $\Omega$.}
		\end{equation}
	\end{lemma}
	\begin{proof}
		Let us begin by choosing $ T_{a,b}(u_n) \varphi \in BV(\Omega) \cap L^\infty(\Omega) $, with $ 0 < a < b < \infty $ and $ 0 \leq \varphi \in C^1_c(\Omega) $, as a test function in \eqref{exp f test n}. Then we have \begin{equation}\label{prob con t e fi}
			\into \left(z_n, DT_{a,b}(u_n)\right) \varphi + \into z_n \cdot \nabla \varphi T_{a,b}(u_n) = \into f_n T_{a,b}(u_n)\varphi.
		\end{equation}
		Our goal is letting $n$ go to infinity in \eqref{prob con t e fi}. For the first integral on the left hand side using \eqref{sci grad} and \eqref{conv TAB}, we gain
		\begin{multline}\label{sci 1int}
			\frac{1}{m+1} \into \left|DT_{a,b}(u)^{m+1}\right| \varphi\\
			\le \liminf_{n \to \infty} \frac{1}{m+1}\into \left|DT_{a,b}(u_n)^{m+1}\right| \varphi\stackrel{\eqref{sol 2 n}}{=} \liminf_{n \to \infty} \into \left(z_n, DT_{a,b}(u_n)\right)\varphi.
		\end{multline}
		For the second integral on the left hand side of \eqref{prob con t e fi}, we use \eqref{con deb zn to z} and \eqref{conv TAB}, we can assert
		\begin{equation}\label{conv 2int}
			\into z \cdot \nabla \varphi T_{a,b}(u) = \lim_{n \to \infty} \into z_n \cdot \nabla \varphi T_{a,b}(u_n).
		\end{equation}
		Finally, for the right hand side we apply Lebesgue Theorem and we prove
		\begin{equation}\label{conv 3int}
			\into fT_{a,b}(u)\varphi = \lim_{n \to \infty} \into f_n T_{a,b}(u_n)\varphi.
		\end{equation}
		Putting together \eqref{sci 1int}, \eqref{conv 2int} and \eqref{conv 3int}, it yields
		\begin{equation*}
			\frac{1}{m+1} \into |DT_{a,b}(u)^{m+1}|\varphi +\into z \cdot \nabla \varphi T_{a,b}(u) \le \into fT_{a,b}(u)\varphi.
		\end{equation*}
		Recalling that $T_{a,b}(u) \in \mathcal{X}^z(\Omega)$, we can affirm
		\begin{equation*}
			\frac{1}{m+1} \into |DT_{a,b}(u)^{m+1}|\varphi \le \into \left(z, DT_{a,b}(u)\right)\varphi,
		\end{equation*}
		where $\left(z, DT_{a,b}(u)\right)$ is a well defined distribution as stated in \eqref{pair l1}.
		\\We observe that $ \left(z, DT_{a,b}(u)\right) $ is a nonnegative distribution and, therefore, a measure. Note that we can assure that $T_{a,b}(u) \in BV^z(\Omega)$. Consequently, the proof is complete.
	\end{proof}
	
	\subsection{Function $u$ belongs to $DTBV(\Omega)$}
	This subsection is devoted to check that  $\mathcal{H}^{N-1}(J_u)=0$. The proof of this fact is complicated since $z$ is no longer a bounded vector filed. Thus, we have to work with the approximate fields $z_n$.
	
	We begin by proving that the pairing $(z_n, DT_b(u_n)^{m+1})$ defines a nonnegative measure for every $n \in \mathbb{N}$ and for a.e. $b > 0$.
	\begin{lemma}\label{lem zntbunm+1 nonne}
		Assume $m>0$ and $0 \le f \in L^1(\Omega)$. For each $n\in\N$, let us consider a solution $u_n$ to  problem \eqref{Prob trans tronc} with vector field $w_n$. Then, for a.e. $b>0$, \begin{equation}\label{pair zv nonne}
			\left(z_n, DT_b(u_n)^{m+1}\right) \quad \text{is a nonnegative measure.}
		\end{equation}
	\end{lemma}
	\begin{proof}
		Let us take $0<a<b<\infty$ and $T_{a,b}(u_n)^{m+1} \varphi\in BV(\Omega) \cap L^\infty(\Omega)$ as a test function in \eqref{exp f test n}, with $0\le \varphi \in C^1_c(\Omega)$,  we have
		\begin{equation}\label{= per la pos di znvn}
			\int_\Omega \left(z_n, DT_{a,b}(u_n)^{m+1}\right) \varphi + \int_\Omega z_n \cdot \nabla \varphi T_{a,b}(u_n)^{m+1}= \int_\Omega f_n T_{a,b}(u_n)^{m+1} \varphi,
		\end{equation}
		so to achieve our goal, we need to take limits as $a$ tends to $0$ in the previous equality.
		\\We may write the first integral on left hand side of \eqref{= per la pos di znvn}, as
		\begin{equation}\label{= pair T2m+1}
			\begin{aligned}
				\left(z_n, DT_{a,b}(u_n)^{m+1}\right)\stackrel{\eqref{pair con der teta}}{=}& \theta \left(z_n, DT_{a,b}(u_n)^{m+1},x\right) \left|DT_{a,b}(u_n)^{m+1}\right| \\ \stackrel{\eqref{ch wj}}{=} &(m+1)T_{a,b}(u_n)^{m} \, \theta\left(z_n, DT_{a,b}(u_n)^{m+1}, x\right) |DT_{a,b}(u_n)|\\ = \ &(m+1)T_{a,b}(u_n)^{m} \, \theta\left(z_n, DT_{a,b}(u_n), x\right) |DT_{a,b}(u_n)|  \\
				\stackrel{\eqref{pair con der teta}}{=}  & (m+1)T_{a,b}(u_n)^m \left(z_n, DT_{a,b}(u_n)\right) \\ \stackrel{\eqref{sol 2 n}}{=}& T_{a,b}(u_n)^m\left|DT_{a,b}(u_n)^{m+1}\right| \\ \stackrel{\eqref{ch wj}}{=}  &\frac{m+1}{2m+1}\left|DT_{a,b}(u_n)^{2m+1}\right|.
			\end{aligned}
		\end{equation}
		where in third line we used Lemma \ref{prop 4.5 cdc}.
		\\We point out that the function $ s \mapsto s^{\frac{2m+1}{m+1}} $ is Lipschitz and that $ u_n^{m+1} \in BV(\Omega) $, we apply Lemma \ref{lem c r} to deduce that $ T_b(u_n)^{2m+1} \in BV(\Omega) \cap L^\infty(\Omega) $. Taking the limit as $a$ tends to $0$ and using \eqref{sci grad}, we obtain
		\begin{equation}\label{per la pos}
			\begin{aligned}
				\frac{m+1}{2m+1} \int_\Omega \left|DT_b(u_n)^{2m+1}\right| \varphi \le &\liminf_{a \to 0} \frac{m+1}{2m+1} \int_\Omega \left| DT_{a,b}(u_n)^{2m+1}\right|  \varphi \\ \stackrel{\eqref{= pair T2m+1}}{=} &\liminf_{a \to 0} \int_\Omega \left(z_n, DT_{a,b}(u_n)^{m+1}\right) \varphi.
			\end{aligned}
		\end{equation}
		For the second integral on left hand side and the integral on right hand side of \eqref{= per la pos di znvn}, we use Lebesgue's Theorem. Gathering these limits with  \eqref{per la pos}, we arrive at
		\begin{equation*}
			\begin{aligned}
				\frac{m+1}{2m+1} \int_\Omega \left|DT_b(u_n)^{2m+1}\right| \varphi \le - &\int_\Omega z_n \cdot \nabla \varphi T_b(u_n)^{m+1} + \int_\Omega f_n T_b(u_n)^{m+1} \varphi \\ \stackrel{\eqref{def pair}}{=} &\int_\Omega \left(z_n, D T_b(u_n)^{m+1}\right) \varphi,
			\end{aligned}
		\end{equation*}
		which implies \eqref{pair zv nonne} and we conclude the proof.
	\end{proof}
	Now we present some properties of the measure $ \left(z, DT_{a,b}(u)\right) $ and the generalized pairing $ \left(z, D\chi_{\{u > k\}}\right) $ with $k>0$.
	\begin{lemma}\label{lemma delle misure}
		Assume $m>0$ and $0 \le f \in L^1(\Omega)$. Let us consider $u_n$ a                  solution to the problem \eqref{Prob trans tronc} with vector field $w_n$. Then for a.e. $k>0$, we have
		\begin{equation}\label{zdchik e mis}
			(z, D\chi_{\{u>k\}}) \in \mathcal{M}(\Omega),
		\end{equation}
this measure is nonnegative and satisfies
		\begin{equation*} 
			\int_\Omega (z, D\chi_{\{u>k\}})\le \int_{\{u>k\}}f,
		\end{equation*}
		Moreover,
		\begin{equation}\label{conv mis}
			\lim_{n\to\infty}\int_\Omega (z_n, D\chi_{\{u_n>k\}})\varphi
			=\int_\Omega (z, D\chi_{\{u>k\}}) \varphi,
		\end{equation}
		for every $\varphi\in C_c^1(\Omega)$.
		
		In addition, for a.e. $0<a<b< \infty$ and $k>0$, it holds
		\begin{equation}\label{campi troncati}
			z_{a,b}, z_k \in \DM(\Omega),
		\end{equation}
		where $z_{a,b}:=z\chi_{\{a < u < b\}}$ and $z_k:= z \chi_{\{u \le k\}}$.
	\end{lemma}
	\begin{proof}
		We initially observe that, for a.e. $0 < a < b < \infty$ and for a.e. $k > 0$,
		\begin{equation*} 
			\chi_{\{a < u < b\}}, \chi_{\{u \le k\}}, \chi_{\{u > k\}} \in BV(\Omega),
		\end{equation*}
		as established in subsection \ref{subs BV}.
		\medskip
		\\The next step is to show \eqref{zdchik e mis}. We recall that $z \in \mathcal{DM}^1(\Omega)$ and $\operatorname{div} z\in L^1(\Omega)$; therefore, thanks to Proposition \ref{prop z e div in l1}, it is readily seen that $\chi_{\{u>k\}} \in \mathcal{X}^z(\Omega)$. As a consequence, the pairing $\left(z, D\chi_{\{u>k\}}\right)$ is a well-defined distribution for a.e. $k>0$.
		\\We now proceed to prove that $\chi_{\{u>k\}} \in BV^z(\Omega)$, for a.e. $k>0$. To demonstrate this property, we adopt a strategy inspired by the proof in \cite[Theorem 4.1]{MST}. Let us fix $b>0$, consider $ 0 < j < b^{m+1} $ and set $ v_n := T_b(u_n)^{m+1} $. Next take $ h_{j, \varepsilon}(v_n) \varphi \in BV(\Omega) \cap L^\infty(\Omega) $ as a test function in \eqref{exp f test n}, where $ h_{j,\varepsilon} $ is the function defined in \eqref{hkeps} and $ 0 \leq \varphi \in C^1_c(\Omega) $, it gives
		\begin{equation*}
			\int_\Omega \left(z_n, D h_{j, \varepsilon}(v_n)\right) \varphi =-\int_\Omega z_n \cdot \nabla \, \varphi h_{j,\varepsilon}(v_n) + \int_\Omega f_n h_{j,\varepsilon}(v_n)\varphi.
		\end{equation*}	
		Firstly, we take the limit as $n$ tends to infinity, for the first integral on the right hand side we use \eqref{con deb zn to z} and \eqref{conv TAB}, while for the second integral we use Lebesgue's Theorem, so we get
		\begin{equation*}
			\lim_{n \to \infty} 	\int_\Omega \left(z_n, D h_{j, \varepsilon}(v_n)\right) \varphi = - \int_\Omega z \cdot \nabla \varphi h_{j, \varepsilon}(v)+ \int_\Omega fh_{j,\varepsilon}(v) \varphi.
		\end{equation*}
		As $\varepsilon$ tends to $0$, we get
		\begin{equation*}
			h_{j,\varepsilon}(v) \to \chi_{\{v\le j\}} \quad \text{a.e. $x \in \Omega$,}
		\end{equation*}
		for almost all $j\in(0,b^{m+1})$.
		Putting $k:=j^{\frac{1}{m+1}}$, it implies
		\begin{equation*}
			h_{j,\varepsilon}(v) \to \chi_{\{u\le k\}} \quad \text{a.e. $x \in \Omega$,}
		\end{equation*}
		for almost all $k\in (0,b)$ and, due to the arbitrariness of $b$, for almost all $k>0$.
		Employing Lebesgue's Theorem, we assert
		\begin{equation*}
			\lim_{\varepsilon \to 0} \lim_{n \to \infty} 	\int_\Omega \left(z_n, D h_{j, \varepsilon}(v_n)\right) \varphi = - \int_\Omega z \cdot \nabla \varphi \chi_{\{u\le k\}} + \int_\Omega f \chi_{\{u\le k\}} \varphi,
		\end{equation*}
		and, having in mind \eqref{sol 1} and \eqref{pair l1}, we gain
		\begin{equation}\label{per la misura}
			-	\lim_{\varepsilon \to 0} 	\lim_{n \to \infty} 	\int_\Omega \left(z_n, D h_{j, \varepsilon}(v_n)\right) \varphi = \int_\Omega \left(z, D\chi_{\{u>k\}}\right) \varphi.
		\end{equation}
		We point out that the above arguments also leads to
		\begin{equation}\label{per la misura-n}
			-	\lim_{\varepsilon \to 0} 	\int_\Omega \left(z_n, D h_{j, \varepsilon}(v_n)\right) \varphi = \int_\Omega \left(z_n, D\chi_{\{u_n>k\}}\right) \varphi.
		\end{equation}
		
		\medskip
		
		Now we choose $T_\varepsilon\left(v_n -T_j(v_n)\right)$ as test function in \eqref{exp f test n}, obtaining
		\begin{equation}\label{trunc test n}
		\begin{aligned}
				&\int_\Omega \left(z_n, DT_\varepsilon\left(v_n -T_j(v_n)\right) \right)\\
			=  &\int_{\partial \Omega} T_\varepsilon\left(v_n-T_j(v_n)\right) [z_n,\nu]\, d\mathcal H^{N-1} + \int_\Omega f_n T_\varepsilon\left(v_n-T_j(v_n)\right) \\
			\le &-\int_{\partial \Omega} T_\varepsilon\left(v_n-T_j(v_n)\right)  \left(u_n^\Omega\right)^m d\mathcal H^{N-1} + \varepsilon \int_{\{v_n>j\}}f\\
			\le &\varepsilon \int_{\{v_n>j\}}f
		\end{aligned}
		\end{equation}
		dropping a nonpositive term.
		
		On the other hand, from \cite[Lemma 5.3]{BOPS}, we notice that
		\begin{equation}\label{pair zv con der}
			\begin{aligned}
				-\left(z_n, Dh_{j,\varepsilon}(v_n)\right) \stackrel{\eqref{pair con der teta}}{=} &\theta\left(z_n, D\left(-h_{j,\varepsilon}(v_n)\right), x \right) \left|Dh_{j,\varepsilon}(v_n)\right| \\ \stackrel{\eqref{ch wj}}{=}-&h_{j,\varepsilon}'(v_n) \, \theta(z_n, Dv_n, x)\left|Dv_n\right|\\ \stackrel{\eqref{pair con der teta}}{=} -&h_{j,\varepsilon}'(v_n) \left(z_n, Dv_n\right),
			\end{aligned}
		\end{equation}
		where in the second equality we used Lemma \ref{prop 4.5 cdc}.
		\medskip
		Hence, we deduce that
		\begin{equation}\label{dis per prov pair mis}
			\begin{aligned}
				\left|\int_\Omega \left(z_n, D h_{j, \varepsilon}(v_n)\right) \varphi\right| \stackrel{\eqref{pair zv con der}}{=} &\left| \int_\Omega h_{j, \varepsilon}'(v_n) (z_n, Dv_n) \varphi\right| \\ \stackrel{\eqref{pair zv nonne}}{\le} &\|\varphi\|_{L^\infty(\Omega)} \int_\Omega \left|h_{j,\varepsilon}'(v_n) \right| (z_n,Dv_n)\\ \stackrel{\eqref{ch wj}}{=} &\frac{\|\varphi\|_{L^\infty(\Omega)}}{\varepsilon} \int_\Omega \left(z_n, DT_\varepsilon\left(v_n -T_j(v_n)\right) \right) \\ \stackrel{\eqref{trunc test n}}{\le} & \|\varphi\|_{L^\infty(\Omega)} \int_{\{v_n>j\}}f.
			\end{aligned}
		\end{equation}
		Several consequences can be inferring from this inequality. First, applying \eqref{per la misura},  \eqref{dis per prov pair mis} we get
		\begin{equation}\label{desig}
			\left|\int_\Omega \left(z, D\chi_{\{u>k\}}\right) \varphi \right|\le \|\varphi\|_{L^\infty(\Omega)} \int_{\{u>k\}}f\le  \|\varphi\|_{L^\infty(\Omega)}  \|f\|_{L^1(\Omega)},
		\end{equation}
		which means that $(z, D\chi_{\{u>k\}})$ is a distribution of rank zero, i.e. it is a Radon measure for a.e. $k>0$.
		Therefore, $\chi_{\{u>k\}}, \chi_{\{u \le k\}} \in BV^z(\Omega)$ and from Lemma \ref{lemma per mis l1}, it holds
		\begin{equation}\label{Leib con divzk}
			\left(z, D\chi_{\{u \le k\}}\right)=-\chi_{\{u \le k\}}\operatorname{div}z+\operatorname{div}z_k \stackrel{\eqref{sol 1}}{=} \chi_{\{u \le k\}}f +\operatorname{div}z_k \quad \text{as measures in $\Omega$,}
		\end{equation}
		hence $z_k \in \DM(\Omega)$. Since $z_k \in \DM (\Omega)$, it follows that $\operatorname{div} z_k \ll \mathcal{H}^{N-1}$. Consequently, we deduce that the measure
		\begin{equation}\label{z><k HN-1}
			\left(z, D\chi_{\{u>k\}}\right)=-\left(z, D\chi_{\{u\le k\}}\right) \ll \mathcal{H}^{N-1}.
		\end{equation}
Moreover, the measure $\left(z, D \chi_{\{u>k\}} \right)$ is nonnegative for a.e. $k>0$ since it arises as the limit of the measures $-h'_{j, \varepsilon}(v_n)(z_n, Dv_n)$, which are themselves nonnegative, as demonstrated in Lemma \ref{lem zntbunm+1 nonne}.
		
		\medskip
		Going back to \eqref{desig}, we also have
		\begin{multline}\label{mi servee}
			\int_\Omega \left(z, D\chi_{\{u>k\}}\right)\\
			= \sup\left\{\into \left(z, D \chi_{\{u > k\}}\right)\varphi \, : \, \varphi \in C^1_c(\Omega), 0\le \varphi \le 1 \right\}
			\le  \int_{\{u>k\}}f,
		\end{multline}
		\medskip
		On account of \eqref{per la misura-n}, it also follows from \eqref{dis per prov pair mis} that
		\[\left|\int_\Omega \left(z_n, D\chi_{\{u_n>k\}}\right)\varphi\right|\le \|\varphi\|_{L^\infty(\Omega)}\int_{\{u>k\}}f,\]
		for all $\varphi\in C_c^1(\Omega)$. Hence, each $\left(z_n, D\chi_{\{u_n>k\}}\right)$ is a nonnegative measure and
		\begin{equation*} 
			\int_\Omega \left(z_n, D\chi_{\{u_n>k\}}\right)\le \int_{\{u>k\}}f\le \|f\|_{L^1(\Omega)} \,,
		\end{equation*}
		holds for every $n\in\N$.
		\medskip
		\\The next step is to demonstrate that $z_{a,b} \in \DM(\Omega)$ for almost every $0<a<b<\infty$. Fix $0<a<b<k$. From the previous step, we know that $z_k \in \DM(\Omega)$, furthermore, based on the discussion in Subsection \ref{subs BV}, we also have $\chi_{\{a<u<b\}}\in BV(\Omega)\cap L^\infty(\Omega)$, hence applying Lemma \ref{lem uz in DM} then directly yields \eqref{campi troncati}, which concludes the proof.
		
		\medskip
		
		It remains to check \eqref{conv mis}. This fact is straightforward from the definition of pairing and the convergences $f_n\to f$ strongly in $L^1(\Omega)$ and $\chi_{\{u_n>k\}}\to \chi_{\{u>k\}}$ strongly in any $L^q(\Omega)$ for every $1 \le q < \infty$ and for almost all $k>0$.
	\end{proof}
	\begin{remark}\label{rem nonne e con zdchi}
		For the measures $\left(z, D\chi_{\{u>k\}}\right)$, we deduce from \eqref{mi servee} that
		\begin{equation}\label{conv forte di zchiu>j}
			\lim_{k \to \infty} \int_\Omega \left(z, D\chi_{\{u>k\}}\right) =0,
		\end{equation}
		from where we easily obtain
		\begin{equation*} 
			\lim_{k \to \infty} \int_\Omega \left(z, D\chi_{\{u>k\}}\right) \varphi =0 \quad \text{ for every } \varphi \in C^1(\overline\Omega)\,.
		\end{equation*}
		To see \eqref{conv forte di zchiu>j}, just take the limit as $k$ tends to infinity in \eqref{mi servee}, this limit exists because $f \in L^1(\Omega)$ and $u$ is a.e. finite, due to \eqref{u fqo}.
	\end{remark}
	
	The following Lemma establishes a fundamental inequality which plays a key role in proving that our solution $u$ belongs to the space $DTBV(\Omega)$.
	\begin{lemma} 
		Assume $m>0$ and $0 \le f \in L^1(\Omega)$. Let us consider $u_n$ a solution to the problem \eqref{Prob trans tronc} with associated vector field $w_n$. Then, for a.e. $0<a<b<\infty$, for a.e. $k>b$ it holds
		\begin{equation}\label{le sol 2 con zk}
			\frac{1}{m+1}\left|DT_{a,b}(u)^{m+1}\right| \le \left(z_k, DT_{a,b}(u)\right) + \left(b-T_{a,b}^*(u)\right)\left(z, D\chi_{\{u>k\}}\right),
		\end{equation}
		as measures in $\Omega$.
	\end{lemma}
	\begin{proof}
		Recall that in Subsection \ref{subs aprx p}, we stated that for every $n \in \mathbb{N}$, the vector field $z_n = u_n^m w_n$ belongs to $\DM(\Omega)$. Moreover, for a.e. $k > 0$, the function $\chi_{\{u_n \le k\}} \in BV(\Omega) \cap L^\infty(\Omega)$. Consequently, Lemma \ref{lem uz in DM} allows us to deduce that $z_n \chi_{\{u_n \le k\}} \in \DM(\Omega)$. Then, using the Definition of the pairing operator (see \eqref{def pair}), and choosing any $0 \le \varphi \in C^1_c(\Omega)$, we obtain
		\begin{equation}\label{n per le sol 2 con zk}
			\into \left(z_n \chi_{\{u_n \le k\}}, DT_{a,b}(u_n)\right) \varphi + \int_{\{u_n \le k\}} T_{a,b}(u_n) z_n \cdot \nabla \varphi = - \into \operatorname{div}\left(\displaystyle z_n \chi_{\{u_n \le k\}}\right) T_{a,b}(u_n) \varphi.
		\end{equation}
		Our goal is to pass to the limit as $n \to \infty$ in the previous identity. We begin by observing that, since $z_n \chi_{\{u_n \le k\}} \in \DM(\Omega)$ and $T_{a,b}(u_n) \in DTBV(\Omega)$, as a consequence of \cite[Lemma 5.3]{GMP}, we have
		$$
		\left(z_n \chi_{\{u_n \le k\}}, DT_{a,b}(u_n) \right) = \left(z_n \chi_{\{u_n \le k\}} \chi_{\{a < u_n < b\}}, DT_{a,b}(u_n)\right) \quad \text{as measures in } \Omega.
		$$
		It follows from the inequalities $k > b > a$ that $z_n \chi_{\{u_n \le k\}} \chi_{\{a < u_n < b\}} = z_n \chi_{\{a < u_n < b\}}$, and thus  the previous identity yields
		\begin{equation}\label{unosss}
			\begin{aligned}
				\frac{1}{m+1}\into \left|DT_{a,b}(u)^{m+1}\right|\varphi \stackrel{\eqref{sci grad}}{\le} & \liminf_{n \to \infty}\frac{1}{m+1}\into \left|DT_{a,b}(u_n)^{m+1} \right|\varphi \\ \stackrel{\eqref{sol 2 n}}{=}&\liminf_{n \to \infty}\into \left(z_n \chi_{\{a < u_n < b\}}, DT_{a,b}(u_n) \right)\varphi \\= &\liminf_{n \to \infty}	\into \left(z_n \chi_{\{u_n \le k\}}, DT_{a,b}(u_n) \right)\varphi.
			\end{aligned}
		\end{equation}
		\\As for the second integral on the left hand side of \eqref{n per le sol 2 con zk}, we know by Corollary \ref{cor ext u} that $u_n^m \to u^m$ strongly in $L^1(\Omega)$, therefore, for almost every $k > 0$, we have $\chi_{\{u_n \le k\}} \to \chi_{\{u \le k\}}$ strongly in $L^q(\Omega)$, for every $1\le q <\infty$. Thus, taking into account \eqref{con deb zn to z}, we obtain
		\begin{equation}\label{doooos}
			\int_{\{u \le k\}} T_{a,b}(u) z \cdot \nabla \varphi = \lim_{n \to \infty} \int_{\{ u_n \le k\}} T_{a,b}(u_n) z_n \cdot \nabla \varphi.
		\end{equation}
		Concerning the right hand side of \eqref{n per le sol 2 con zk}, we get
		\begin{equation*}
			\begin{aligned}
					- &\lim_{n \to \infty} \into \operatorname{div}\left(\displaystyle z_n \chi_{\{u_n \le k\}}\right) T_{a,b}(u_n)\varphi\\
				\stackrel{\eqref{leib rule dm infty}}{=}
				&\lim_{n \to \infty}-\int_{\{u_n \le k\}} \operatorname{div}z_n T_{a,b}(u_n) \varphi -\into \left(z_n, D\chi_{\{u_n \le k\}}\right) T_{a,b}(u_n)\varphi\\
				\stackrel{\eqref{sol 1 n}}{=}
				&\lim_{n \to \infty}\int_{\{u_n \le k\}} f_n T_{a,b}(u_n) \varphi -\into \left(z_n, D\chi_{\{u_n \le k\}}\right) T_{a,b}(u_n)\varphi.
			\end{aligned}
		\end{equation*}
		We recall that $\left(z_n, D\chi_{\{u_n \le k\}}\right)=-\left(z_n, D\chi_{\{u_n > k\}}\right)$ is absolutely continuous with respect to $\left|D \chi_{\{u_n \le k\}}\right|$ (by \eqref{pair bel mis}) which is concentrated on the set $\{u_n=k\}$, so that
		\begin{equation}\label{quasi secondo int per le sol 2 zk}
	\begin{aligned}
				-&\lim_{n \to \infty} \into \operatorname{div}\left(\displaystyle z_n \chi_{\{u_n \le k\}}\right) T_{a,b}(u_n)\varphi\\
		=&\lim_{n \to \infty}\into f_n \chi_{\{u_n \le k\}} T_{a,b}(u_n) \varphi +b \into \left(z_n, D\chi_{\{u_n > k\}}\right)\varphi\\
		=&\into f \chi_{\{u \le k\}} T_{a,b}(u) \varphi + b \into \left(z,D\chi_{\{u>k\}}\right)\varphi,
	\end{aligned}
		\end{equation}
		as a consequence of Lebesgue's Theorem and Lemma \ref{lemma delle misure}.
		\\Furthermore, using \eqref{Leib con divzk}, equation \eqref{quasi secondo int per le sol 2 zk} implies
		\begin{equation}\label{terzo int per le sol 2 zk}
			\begin{aligned}
				-& \lim_{n \to \infty} \into \operatorname{div}\left(\displaystyle z_n \chi_{\{u_n \le k\}}\right) T_{a,b}(u_n)\varphi\\
				= & \into f \chi_{\{ u \le k\}}T_{a,b}(u)\varphi + \into b \left(z, D\chi_{\{u>k\}}\right) \varphi\\
				=& \into f \chi_{\{u \le k\}}T_{a,b}(u)\varphi -\into T_{a,b}^*(u) \left(z,D\chi_{\{u\le k\}}\right)\varphi+ \into \left(b-T_{a,b}^*(u)\right)\left(z, D\chi_{\{u>k\}}\right)\varphi\\
				\stackrel{\eqref{Leib con divzk}}{=}&
				-\into \operatorname{div}z_k T_{a,b}^*(u)\varphi+\into \left(b-T_{a,b}^*(u)\right) \left(z,D\chi_{\{u>k\}}\right)\varphi.
			\end{aligned}
		\end{equation}
		Combining \eqref{unosss}, \eqref{doooos} and \eqref{terzo int per le sol 2 zk}, the identity \eqref{n per le sol 2 con zk} becomes
		\begin{equation*}
			\begin{aligned}
				\frac{1}{m+1} &\into \left|D T_{a,b}(u)^{m+1}\right|\varphi + \int_{\{u \le k\}} T_{a,b}(u) z \cdot \nabla \varphi \\ \le - &\into \operatorname{div}z_k T_{a,b}^*(u)\varphi + \into \left(b-T_{a,b}^*(u)\right) \left(z,D\chi_{\{u>k\}}\right) \varphi,
			\end{aligned}
		\end{equation*}
		recalling the Definition of pairing by Anzellotti \eqref{def pair}, we derive \eqref{le sol 2 con zk} and this concludes the proof.
	\end{proof}

	Our next result is the main one in this subsection since it shows that $\mathcal{H}^{N-1}(J_u)=0$.
	\begin{lemma}\label{lem u j v}
		Let us consider $u$ the function provided in Corollary \ref{cor ext u}. Then, for a.e. $0<a<b< \infty$, it holds \begin{equation*}
			\mathcal{H}^{N-1}\left(J_{T_{a,b}(u)}\right)=0.
		\end{equation*}
		As a consequence, $\mathcal{H}^{N-1}(J_u)=0$.
	\end{lemma}
	\begin{proof}
		The proof is an adaptation of \cite[Lemma 5.9]{GMP} to our setting.
		\\ Let us fix $0<a<b<\infty$, from \cite[Proposition 3.69]{AFP} we have
		$$J_{T_{a,b}(u)}=J_{T_{a,b}(u)^{m+1}},$$
		owing to $m>0$. Moreover, on this set the corresponding orientations coincide: $\nu_{T_{a,b}(u)}=\nu_{T_{a,b}(u)^{m+1}}$.
		\\Recalling that $J_{T_{a,b}(u)}$ is a countably $\mathcal{H}^{N-1}$-rectifiable set, there exist regular hypersurfaces $\xi_i$ ($i\in I$, $I$ countable) such that $$\mathcal{H}^{N-1}\left(J_{T_{a,b}(u)} \setminus\bigcup_{i \in I} \xi_i \right)=0.$$
		In order to show $\mathcal{H}^{N-1}(J_{T_{a,b}(u)})=0$, we fix $i\in I$ and see that $T_{a,b}(u)^+=T_{a,b}(u)^-$ on $\xi_i$. We stress that we may assume that $\xi_i\Subset\Omega$ since if we prove $T_{a,b}(u)^+=T_{a,b}(u)^-$ on each set $\xi_i\cap\{x\in\Omega\,:\, dist(x,\partial\Omega)>\epsilon\}$, then we are done. Thus, we may consider a regular open set $\omega_i \Subset \Omega$ such that $\xi_i\subset \partial \omega_i$.
		For a.e. $k>b$, using $f\chi_{\{u \le k\}} \in L^1(\Omega)$, \eqref{Leib con divzk} and \eqref{t div su ip}, we get
		\begin{equation}\label{due volte di divzk}
			\begin{aligned}
				\left(z, D \chi_{\{u>k\}}\right) \res \xi_i= &\left(f\chi_{\{u \le k\}}+ \left(z, D \chi_{\{u>k\}}\right) \right) \res \xi_i \\ = &-\operatorname{div}z_k\res \xi_i \\ = &\left(-[z_k,\nu_{\xi_i}]^+ + [z_k, \nu_{\xi_i}]^- \right)\mathcal{H}^{N-1} \res \xi_i,
			\end{aligned}
		\end{equation}
		which we write
		\begin{equation}\label{primero}
			[z_k,\nu_{\xi_i}]^- \mathcal{H}^{N-1}\res{\xi_i}= [z_k,\nu_{\xi_i}]^+ \mathcal{H}^{N-1}\res{\xi_i} + \left(z, D\chi_{\{u>k\}}\right)\res{\xi_i}.
		\end{equation}
		\\On the other hand, we perform the following manipulations
		\begin{equation}\label{per salto nullo}
			\begin{aligned}
				&	\frac{1}{m+1} \left|DT_{a,b}(u)^{m+1}\right|\res \xi_i\\
				\stackrel{\eqref{le sol 2 con zk}}{\le} & \left(\left(z_k, DT_{a,b}(u)\right) + \left(b-T_{a,b}^*(u)\right)\left(z, D\chi_{\{u>k\}} \right) \right)\res \xi_i \\ \stackrel{\eqref{leib rule dm infty}}{=}& \left(-T_{a,b}(u)^* \operatorname{div}z_k + \operatorname{div}\left(\displaystyle T_{a,b}(u)z_k\right)+ \left(b-T_{a,b}^*(u)\right)\left(z, D\chi_{\{u>k\}} \right) \right) \res \xi_i \\ \stackrel{\eqref{t div su ip}}{=}& \left(-T_{a,b}^*(u) [z_k,\nu_{\xi_i}]^+ + T_{a,b}^*(u)[z_k, \nu_{\xi_i}]^- + [T_{a,b}(u)z_k,\nu_{\xi_i}]^+ -[T_{a,b}(u)z_k, \nu_{\xi_i}]^- \right) \mathcal{H}^{N-1}\res{\xi_i}\\  + &\left(b-T_{a,b}^*(u)\right)\left(z, D\chi_{\{u>k\}} \right) \res \xi_i\\ \stackrel{\eqref{saltando sul bordo}}{=} &\left(-T_{a,b}(u)^*[z_k,\nu_{\xi_i}]^+ +T_{a,b}^*(u) [z_k, \nu_{\xi_i}]^- + T_{a,b}(u)^+[z_k,\nu_{\xi_i}]^+ - T_{a,b}(u)^-[z_k, \nu_{\xi_i}]^- \right) \mathcal{H}^{N-1}\res{\xi_i} \\ + &\left(b-T_{a,b}^*(u)\right)\left(z, D\chi_{\{u>k\}} \right) \res \xi_i.
			\end{aligned}
		\end{equation}
		Hence, on the hyperplane $\xi_i$, having in mind \eqref{due volte di divzk}, the inequality \eqref{per salto nullo} becomes
		\begin{equation}\label{nuovo cal senza j}
			\begin{aligned}
				& \frac{1}{m+1} \left|DT_{a,b}(u)^{m+1} \right|\res \xi_i\\
				\le &(b-T_{a,b}^-(u)) \left(z, D\chi_{\{u>k\}}\right)\res \xi_i +\left(T_{a,b}^+(u)-T_{a,b}^-(u)\right)[z_k,\nu_{\xi_i}]^+\mathcal{H}^{N-1}\res{\xi_i}
				\\ =& (b-T_{a,b}^-(u)) \left(z,D\chi_{\{u>k\}} \right)\res \xi_i+ \left(T_{a,b}^+(u)-T_{a,b}^-(u)\right)[z_k,\nu_{\xi_i}]^+\chi_{\{a < u < b\}}^+ \mathcal{H}^{N-1}\res{\xi_i}\\ \stackrel{\eqref{saltando sul bordo}}{=} & (b-T_{a,b}^-(u)) \left(z,D\chi_{\{u>k\}} \right)\res \xi_i+ \left(T_{a,b}^+(u)-T_{a,b}^-(u)\right)[z_k\chi_{\{a < u < b\}},\nu_{\xi_i}]^+ \mathcal{H}^{N-1}\res{\xi_i} \\ =&(b-T_{a,b}^-(u)) \left(z,D\chi_{\{u>k\}} \right)\res \xi_i+ \left(T_{a,b}^+(u)-T_{a,b}^-(u)\right)[z\chi_{\{a < u < b\}},\nu_{\xi_i}]^+ \mathcal{H}^{N-1}\res{\xi_i} \\ \stackrel{\eqref{saltando sul bordo}}{=} &(b-T_{a,b}^-(u)) \left(z,D\chi_{\{u>k\}} \right)\res \xi_i\\
				&+ \left(T_{a,b}^+(u)-T_{a,b}^-(u)\right)\left(T_{a,b}(u)^m\right)^+[w\chi_{\{a < u < b\}},\nu_{\xi_i}]^+ \mathcal{H}^{N-1}\res{\xi_i}.
			\end{aligned}
		\end{equation}
		Our next step is to check that $[w\chi_{\{a < u < b\}},\nu_{\xi_i}]^+\le1$. Observe that, as a consequence of Lemmas \ref{lem uz in DM}, \ref{lem ext z} and \ref{lemma delle misure}, we obtain $w\chi_{\{a < u < b\}} \in \DM(\Omega)$ and $\|w \chi_{\{a < u < b\}}\|_{L^\infty(\Omega)^N}\le 1$. Thus, we deduce from \eqref{nuovo cal senza j} that indeed $[w\chi_{\{a < u < b\}},\nu_{\xi_i}]^+\le1$ and so it results
		\begin{equation*}
			\frac{1}{m+1} \left|DT_{a,b}(u)^{m+1} \right|\res \xi_i \le
			\left( b-T_{a,b}^-(u) \right) \left(z,D\chi_{\{u>k\}} \right)\res \xi_i+ \left(T_{a,b}^+(u)-T_{a,b}^-(u)\right)\left(T_{a,b}(u)^m\right)^+\mathcal{H}^{N-1}\res{\xi_i}.
		\end{equation*}
		\\ The passage to the limit when k goes to infinity is justified by the convergence property \eqref{conv forte di zchiu>j}. Applying Lemma \ref{lem c r}, we then obtain
		\begin{equation}\label{tracci min tabm+}
			\frac{1}{m+1} \left|\left(T_{a,b}(u)^{m+1}\right)^{+} -\left(T_{a,b}(u)^{m+1}\right)^{-}\right| \le \left(T_{a,b}^+(u)-T_{a,b}^-(u)\right)\left(T_{a,b}(u)^m\right)^+,
		\end{equation}
		$\mathcal{H}^{N-1}$-a.e. on $\xi_i$.
		\\We now repeat the same reasoning as in the previous step, this time we write \eqref{primero} as
		$$	[z_k,\nu_{\xi_i}]^+ \mathcal{H}^{N-1}\res \xi_i= [z_k,\nu_{\xi_i}]^- \mathcal{H}^{N-1}\res \xi_i - \left(z, D\chi_{\{u>k\}}\right)\res \xi_i.$$
		We now deduce that $$	\frac{1}{m+1} \left|\left(T_{a,b}(u)^{m+1}\right)^{+} -\left(T_{a,b}(u)^{m+1}\right)^{-}\right| \le \left(T_{a,b}^+(u)-T_{a,b}^-(u)\right)\left(T_{a,b}(u)^m\right)^- \quad \text{$\mathcal{H}^{N-1}$-a.e. on $\xi_i$}.$$
		Combining this, with equation \eqref{tracci min tabm+}, this implies
		\begin{equation*}
			\begin{aligned}
				\frac{1}{m+1} \left| \left(T_{a,b}(u)^{m+1}\right)^{+} - \left(T_{a,b}(u)^{m+1}\right)^{-}\right| \le\left(T_{a,b}(u)^+ - T_{a,b}(u)^-\right) \min \left\{\left(T_{a,b}(u)^m\right)^\pm\right\},
			\end{aligned}
		\end{equation*}
		$\mathcal{H}^{N-1}$-a.e. on $\xi_i$.
		\\We arrive at a contradiction with the Mean Value Theorem, since the function $s \mapsto s^m$, with $m > 0$, is increasing. Hence we obtain
		$$0=\mathcal{H}^{N-1}\left(J_{T_{a,b}(u)}\right)=\mathcal{H}^{N-1}\left(S_{T_{a,b}(u)}\right) \quad \text{for a.e. $0<a<b< \infty$,}$$
		and, by Lemma \ref{salti in TBV}, $\mathcal H^{N-1}(S_u^*)=0$.
		Hence, it follows from $J_u \subseteq J_u^* \subseteq S_u^*$ that $$\mathcal{H}^{N-1}\left(J_u\right)=0,$$ and this concludes the proof.
	\end{proof}

	\subsection{Condition (\ref{sol 2})}
	In Subsection \ref{sub 5.3} we have seen one of the inequalities of condition \eqref{sol 2}. Once we get that $u$ has not jump part, this subsection is devoted to check the opposite inequality.
	
	We first prove an important property of the pairing $\left(z, D T_{a,b}(u)^{\alpha}\right)$ for $\alpha > 0$, which will be instrumental in establishing the forthcoming results.
	\begin{lemma}\label{lem zab=z}
		Assume $m>0$ and $0 \le f \in L^1(\Omega)$. Consider that $u_n$ is a solution to the problem \eqref{Prob trans tronc} and $w_n$ its associated vector field. Then for a.e. $0<a<b< \infty$ and for every $\alpha>0$, it holds
		\begin{equation}\label{= pair z e zab}
			\left(z, DT_{a,b}(u)^\alpha\right)=\left(z_{a,b}, DT_{a,b}(u)^\alpha\right) \quad \text{as measures in $\Omega$.}
		\end{equation}
	\end{lemma}
	\begin{proof}
		As illustrated in \cite[Lemma 5.3]{GMP} and \cite[Lemma 2.11]{BOPS}, given $0<a<b<\infty$ and $k>0$, and choosing  $T_{a,b}(u) \in DTBV(\Omega) \cap L^\infty(\Omega)$ and $z_k \in \DM(\Omega)$ (which holds a.e.), we have
		\begin{equation}\label{per = pair zab}
			\left(z_k, DT_{a,b}(u)^\alpha\right)=\left(z_k\chi_{\{a<u<b\}}, DT_{a,b}(u)^\alpha\right) \quad \text{as measures in $\Omega$.}
		\end{equation}
		Now we want let $k$ go to infinity in \eqref{per = pair zab}.
		\\For the pairing on right hand side, we note that for $k>b$, it holds $$z_k\chi_{\{a<u<b\}}=z\chi_{\{a<u<b\}}=z_{a,b}.$$
		Therefore,
		\begin{equation}\label{z_kab=zab}
			\lim_{k \to \infty} \int_\Omega \left(z_k\chi_{\{a<u<b\}}, DT_{a,b}(u)^\alpha\right)\varphi=\int_\Omega \left(z_{a,b}, DT_{a,b}(u)^\alpha\right)\varphi, \quad \text{for all $\varphi \in C_c^1(\Omega)$.}
		\end{equation}
		We turn to deal with the pairing on the left hand side of equation \eqref{per = pair zab}. For every $0 \le \varphi \in C^1_c(\Omega)$, we observe
		\begin{equation}\label{insieme}
			\begin{aligned}
				\int_\Omega
				\left(z_k, DT_{a,b}(u)^\alpha\right)\varphi \stackrel{\eqref{def pair}}{=}  &-\int_\Omega T_{a,b}(u)^\alpha \varphi \operatorname{div}z_k - \int_\Omega  z_k \cdot \nabla \varphi T_{a,b}(u)^\alpha \\ \stackrel{\eqref{ug pair con div l1}}{=} &-\int_\Omega T_{a,b}(u)^\alpha\varphi \chi_{\{u \le k\}} \operatorname{div}z  \\ &- \int_\Omega T_{a,b}(u)^\alpha  \left(z, D\chi_{\{u \le k\}}\right) \varphi-\int_{\{u \le k\}} z \cdot \nabla \varphi T_{a,b}(u)^\alpha \\ \stackrel{\eqref{sol 1}}{=} &\int_{\{u \le k\}}  f T_{a,b}(u)^\alpha \varphi -\int_{\{u \le k\}} z \cdot \nabla \varphi T_{a,b}(u)^\alpha\\
				&+ \int_\Omega T_{a,b}(u)^\alpha \left(z, D\chi_{\{u>k\}}\right) \varphi.
			\end{aligned}
		\end{equation}
		We point out that
		$$\int_\Omega T_{a,b}(u)^\alpha \left(z, D\chi_{\{u \le k\}}\right) \varphi, \int_\Omega T_{a,b}(u)^\alpha \left(z, D\chi_{\{u>k\}}\right) \varphi, $$ are well defined integrals, because $T_{a,b}(u)^\alpha \in BV(\Omega) \cap L^\infty(\Omega)$ and \eqref{z><k HN-1} holds.
		By Remark \ref{rem u fqo} and Lebesgue's Theorem, we obtain
		\begin{equation}\label{primo pezzo}
			\begin{aligned}
				\lim_{k\to \infty} \int_{\{u \le k\}} f T_{a,b}(u)^\alpha\varphi - \int_{\{u \le k\}} z \cdot \nabla \varphi T_{a,b}(u)^\alpha= &\int_\Omega f T_{a,b}(u)^\alpha \varphi - \int_\Omega z \cdot \nabla \varphi T_{a,b}(u)^\alpha \\ \stackrel{\eqref{pair l1}}{=} &\int_\Omega \left(z, DT_{a,b}(u)^\alpha\right)\varphi.
			\end{aligned}
		\end{equation}
		It only remains to prove that
		\begin{equation}\label{pair con chiu<k 0}
			\lim_{k \to \infty} \int_\Omega T_{a,b}(u)^\alpha  \left(z, D\chi_{\{u>k\}}\right) \varphi=0,
		\end{equation}
		and this fact is a consequence of Remark \ref{rem nonne e con zdchi}. Going back to \eqref{insieme} and putting together  \eqref{primo pezzo} and \eqref{pair con chiu<k 0}, it yields
		\[\lim_{k\to\infty}\int_\Omega \left(z_k, DT_{a,b}(u)^\alpha\right)\varphi=\int_\Omega\left(z, DT_{a,b}(u)^\alpha\right)\varphi\,.\]
		This fact and \eqref{z_kab=zab} gives \eqref{= pair z e zab}, and so the proof is complete.
	\end{proof}
	\begin{remark}
		From \eqref{= pair z e zab} for every open set $A \Subset \Omega$ and every $\varphi \in C^1_c(A)$, it gives
		\begin{equation*}
			\left|\langle	\left(z, DT_{a,b}(u)\right), \varphi \rangle \right|= \left|\langle	\left(z_{a,b}, DT_{a,b}(u)\right), \varphi \rangle \right| \stackrel{\eqref{pair bel mis}}{\le}\|\varphi\| _{L^\infty(A)} b^m \int_A |DT_{a,b}(u)|,
		\end{equation*}
		which implies that
		\begin{equation*}
			\left(z, DT_{a,b}(u)\right) \quad \text{is absolutely continuous with respect to $|DT_{a,b}(u)|$.}
		\end{equation*}
	\end{remark}
	Now we are able to show that inequality \eqref{dis per sol 2} is actually an equality.
	\begin{lemma}\label{lem per sol 2}
		Let us consider $u$ the function provided in Corollary \ref{cor ext u} and $z \in L^{\frac{N}{N-1},\infty}(\Omega)^N$ the vector field provided in Lemma \ref{lem ext z}. Then, it holds \eqref{sol 2}.
	\end{lemma}
	\begin{proof}
		It is a straightforward adaptation of \cite[Lemma 5.4]{BOPS}. We begin by considering inequality \eqref{dis per sol 2}:
		\begin{equation}\label{prov u wj}
			\begin{aligned}
				\frac{1}{m+1} \left|DT_{a,b}(u)^{m+1}\right| \le &\left(z, DT_{a,b}(u)\right) \\ \stackrel{\eqref{= pair z e zab}}{=} &\left(z_{a,b}, DT_{a,b}(u)\right) \\ = &\left(T_{a,b}(u)^m w \chi_{\{a < u < b\}}, DT_{a,b}(u)\right) \\ \stackrel{\eqref{uscire sx}}{=} &T_{a,b}(u)^m \left(w \chi_{\{a < u < b\}}, DT_{a,b}(u)\right) \\  \stackrel{\eqref{pair bel mis}}{\le} &T_{a,b}(u)^m \, |DT_{a,b}(u)| \\ \stackrel{\eqref{ch wj}}{=} &\frac{1}{m+1} \left|DT_{a,b}(u)^{m+1}\right| ,
			\end{aligned}
		\end{equation}
		in particular, in the last equality, we used the fact that $u \in DTBV(\Omega)$, as shown in Lemma \ref{lem u j v}. Thus we get \eqref{sol 2} and this concludes the proof.
	\end{proof}
	\begin{remark}
		We emphasize that from \eqref{prov u wj}, one can derive an equivalent formulation of equation \eqref{sol 2}, which is
		\begin{equation*}
			\frac{1}{m+1}\left|DT_{a,b}(u)^{m+1}\right| = \left(z_{a,b}, DT_{a,b}(u)\right) \quad \text{as measures in $\Omega$, for a.e. $0 < a <b < \infty$.}
		\end{equation*}
		Thanks to \eqref{prov u wj}, we are able to provide the following equivalent formulation of \eqref{sol 2}
		\begin{equation*}
			\left(w \chi_{\{a<u<b\}}, DT_{a,b}(u)\right)= \left|DT_{a,b}(u)\right| \quad \text{as measure in $\Omega$, for a.e. $0<a<b<\infty$.}
		\end{equation*}
	\end{remark}
	
	\subsection{Condition on the boundary}
	Before handling the boundary condition, have to establish a regularity result, namely that the trace of $z$ on $\partial \Omega$ is a function in $L^1(\partial \Omega)$.
	\begin{lemma}\label{lem trac z in l1}
		Assume $m>0$ and $0 \le f \in L^1(\Omega)$. Let us consider $u_n$ a solution to problem \eqref{Prob trans tronc}, the associated vector field $w_n$ and $z$ the vector field provided in Lemma \ref{lem ext z}. Then,
		\begin{equation*}
			[z,\nu] \in L^1(\partial \Omega).
		\end{equation*}
		Moreover, the following Gauss-Green formula holds
		\begin{equation}\label{GG form definitiva}
			\into (z, DT_{a,b}(u)^\alpha)- \int_{\partial \Omega} T_{a,b}(u)^\alpha  [z,\nu] \, \mathrm{d} \mathcal{H}^{N-1} = \into fT_{a,b}(u)^\alpha,
		\end{equation}
		for all $\alpha >0$ and a.e. $0<a<b<\infty$; and
		\begin{equation}\label{chi fuori chi dentro}
			[z,\nu]\chi_{\{u \le k \}}= [z_k,\nu] \quad \text{for $\mathcal{H}^{N-1}$-a.e. $ x \in \partial \Omega$, for a.e. $k>0$.}
		\end{equation}
	\end{lemma}
	\begin{proof}
		We start by proving that the vector field $z$ obtained in Lemma \ref{lem ext z} admits a trace in $L^1(\partial \Omega)$. This trace will be obtained as the limit of the traces $[z_n, \nu]$, where we recall that $z_n$ is the vector field founded in \cite[Theorem 3.3]{BOPS}.
		To establish the result, we shall make use of the Dunford-Pettis Theorem (see \cite[Theorem 1.38]{AFP}).
		\\As pointed out in Remark \ref{oss-2}, the normal trace $[z_n,\nu]$ is a nonpositive bounded function. We are showing the uniform bound of the sequence $[z_n,\nu]$ in $\|\cdot\|_{L^1(\partial \Omega)}$. To this end, we apply Green's formula \eqref{GG form} to get
		\begin{equation}\label{stime norm l1 trac zn}
			\|[z_n,\nu]\|_{L^1(\partial \Omega)}=	- \int_{\partial \Omega} [z_n,\nu] \ensuremath{\mathrm d}\mathcal{H}^{N-1}=\into f_n \le \|f\|_{L^1(\Omega)}.
		\end{equation}
		Now we show the equi-integrability condition. To get it, we take $\chi_{\{u_n^m>k^m\}} \in BV(\Omega)\cap L^\infty(\Omega)$ with $k>0$ as a test function in \eqref{exp f test n}, recalling that $(z, D\chi_{\{u>k\}})$ is a nonnegative measure as stated in Remark \ref{rem nonne e con zdchi}, and so we obtain
		\begin{equation*}
			\begin{aligned}
				0 \le \lim_{k \to \infty} \sup_{n \in \mathbb{N}}&\int_{\{u_n^m>k^m\}\cap\partial\Omega}\left(u_n^\Omega\right)^m  \ensuremath{\mathrm d}\mathcal{H}^{N-1} \\ \stackrel{\eqref{sol 3 n}}{\le}\lim_{k \to \infty} \sup_{n \in \mathbb{N}} &\left( \into \left(z_n, D\chi_{\{u_n>k\}}\right)- \int_{\partial \Omega} \chi_{\{u_n^m>k^m\}} [z_n,\nu] \ensuremath{\mathrm d}\mathcal{H}^{N-1} \right)\\ \stackrel{\eqref{GG form}}{=} \lim_{k \to \infty} \sup_{n \in \mathbb{N}}&\into f_n \chi_{\{u_n>k\}}=0,
			\end{aligned}
		\end{equation*}
		 since $0 \le f \in L^1(\Omega)$ and it holds \eqref{u fqo}.
		Hence there exists a nonpositive function $ g \in L^1(\partial \Omega)$ such that \begin{equation}\label{conv debole tracce zn}
			[z_n, \nu] \rightharpoonup g \quad \text{weakly in $L^1(\partial \Omega)$.}
		\end{equation}
		The next step is to prove that $g$ is, in the weak sense, the trace of $z$ on $\partial \Omega$. Let us choose a nonnegative $ \varphi \in C^1({\Omega})\cap C(\overline{\Omega})$ as a test function in \eqref{exp f test n} and let us take limit as $n$ tends to infinity. It yields
		\begin{equation*}
			\int_\Omega z_n \cdot \nabla \varphi - \int_{\partial \Omega} \varphi [z_n,\nu] \, \mathrm{d} \mathcal{H}^{N-1}= \int_\Omega f_n \varphi.
		\end{equation*}
		We apply property \eqref{con deb zn to z} for the first integral, the weak convergence \eqref{conv debole tracce zn} to the second integral and Lebesgue's Theorem to the third integral. It gives
		\begin{equation}\label{GG con z e reg}
			\into z \cdot \nabla \varphi - \int_{\partial \Omega} \varphi g \, \mathrm{d} \mathcal{H}^{N-1}= \into f\varphi.
		\end{equation}
		Thus, we may rewrite $g$ as
		\begin{equation*}
			[z,\nu]:=g,
		\end{equation*}
		so that $[z,\nu]\in L^1(\partial\Omega)$.
		Note that $[z,\nu]$ is nonpositive since it is the weak limit of nonpositive functions. In addition, it follows from \eqref{stime norm l1 trac zn} that
		\begin{equation*} 
		\|[z,\nu]\|_{L^1(\partial \Omega)}=	-\int_\Omega [z, \nu]\, d\mathcal H^{N-1}\le \|f\|_{L^1(\Omega)}.
		\end{equation*}
		
		To establish Gauss-Green formula \eqref{GG form definitiva}, we prove that $[z, \nu]$ is the limit as $k \to \infty$ of $[z_k,\nu]$, the trace of the vector field $z_k$ which belongs to $\DM(\Omega)$ for a.e. $k>0$ (see Lemma \ref{lemma delle misure}). We are showing that $[z_k,\nu]$ is uniformly bounded in $L^1(\partial \Omega)$ with respect to $k>0$.
		For a.e. $k>0$, we can affirm that $\chi_{\{ u_n \le k\}} \in BV(\Omega)\cap L^\infty(\Omega)$, by applying Lemma \ref{lem uz in DM}. So, we deduce that $z_n \chi_{\{u_n \le k\}} \in \DM(\Omega)$ and, as a consequence of Gauss-Green formula \eqref{GG form}, it holds
		\begin{equation}\label{per stime di trac znk}
			\int_{\{ u_n \le k\}} z_n \cdot \nabla \varphi - \int_{\partial \Omega} \varphi [z_n \chi_{\{ u_n \le k\}}, \nu] \, \mathrm{d}\mathcal{H}^{N-1}= -\into \operatorname{div}\left(\displaystyle z_n \chi_{\{ u_n \le k\}}\right) \varphi,
		\end{equation}
		for all $0 \le \varphi \in C^1(\Omega) \cap C(\overline{\Omega})$.
		Our objective is to let $n$ go to infinity in the previous equality.
		\\As for the first integral on the left hand side of \eqref{per stime di trac znk}, we exploit the convergence \eqref{con deb zn to z} and the strong convergence $\chi_{\{ u_n \le k \}} \to \chi_{\{ u \le k \}}$ in any $L^r(\Omega)$, for $1\le r<\infty$. Then we obtain
		\begin{equation}\label{primo int znk}
			\int_{\{ u \le k\}} z \cdot \nabla \varphi = \lim_{n \to \infty} \int_{\{u_n \le k\}}z_n \cdot \nabla \varphi.
		\end{equation}
		Instead, for the second integral on the left hand side of \eqref{per stime di trac znk}, we recall that $[z_n \chi_{\{ u_n \le k \}}, \nu] = \chi_{\{ u_n \le k \}} [z_n, \nu]$, as a consequence of property \eqref{v fuori v dentro}. Therefore, by applying Lemma \ref{lemmino di Sergio}, we deduce
		\begin{equation}\label{secondo int znk}
		\begin{aligned}
				\int_{\partial \Omega} \chi_{\{u \le k\}} \varphi [z,\nu] \, \mathrm{d} \mathcal{H}^{N-1}= &\lim_{n \to \infty} \int_{\partial \Omega} \chi_{\{ u_n \le k\}}\varphi[z_n, \nu] \, \mathrm{d} \mathcal{H}^{N-1}\\
			= &\lim_{n \to \infty} \int_{\partial \Omega} \varphi[z_n \chi_{\{ u_n \le k\}}, \nu] \, \mathrm{d} \mathcal{H}^{N-1}.
		\end{aligned}
		\end{equation}
		Finally, we are analyzing the integral on the right hand side of \eqref{per stime di trac znk}. A Remark is in order. Despite Lemma \ref{lemma delle misure}, we still do not know that
		\begin{equation*}
			\lim_{n\to\infty}\int_\Omega \left(z_n, D\chi_{\{u_n>k\}}\right)\varphi=\int_\Omega \left(z, D\chi_{\{u>k\}}\right)\varphi,
		\end{equation*}
		holds, since $\varphi$ does not vanish on the boundary (see Remark \ref{conv z_n a z} below).
		We compute $\into \operatorname{div}\left(\displaystyle z_n \chi_{\{u_n \le k\}}\right) \varphi$ in another way, we have,
		\begin{equation}\label{terzo int znk}
			\begin{aligned}
				& -\lim_{n \to \infty } \into \operatorname{div}\left(\displaystyle z_n \chi_{\{u_n \le k\}}\right) \varphi\\
				\stackrel{\eqref{leib rule dm infty}, \eqref{sol 1 n}}{=}&
				\lim_{n \to \infty} \into f_n \chi_{\{ u_n \le k\}} \varphi + \into \left(z_n, D \chi_{\{ u_n > k\}}\right) \varphi \\
				=& \into f \chi_{\{ u \le k\}} \varphi + \lim_{n \to \infty} \into \left(z_n, D\chi_{\{ u_n > k\}}\right) \varphi \\
				=&-\int_\Omega \operatorname{div}z_k \varphi -\into \left(z, D \chi_{\{u > k\}}\right) \varphi + \lim_{n \to \infty} \into \left(z_n, D\chi_{\{ u_n > k\}}\right) \varphi.
			\end{aligned}
		\end{equation}
Recalling \eqref{per stime di trac znk}, gathering \eqref{primo int znk}, \eqref{secondo int znk} and \eqref{terzo int znk}, it yields
		\begin{equation*}
			\begin{aligned}
				&\int_{\{ u \le k\}} z \cdot \nabla \varphi  - \int_{\partial \Omega} \chi_{\{u \le k\}} \varphi [z,\nu] \, \mathrm{d} \mathcal{H}^{N-1}\\ = -&\int_\Omega \operatorname{div}z_k \varphi - \into \left(z, D \chi_{\{u > k\}}\right) \varphi+ \lim_{n \to \infty} \into \left(z_n, D\chi_{\{ u_n > k\}}\right) \varphi .
			\end{aligned}
		\end{equation*}
		Appealing again to \eqref{GG form}, it follows that
		\begin{multline}\label{Gallo}
			\int_{\partial \Omega}\varphi[z_k,\nu] \, \mathrm{d}\mathcal{H}^{N-1} \\ =\int_{\partial \Omega} \chi_{\{u \le k\}} \varphi [z,\nu] \, \mathrm{d} \mathcal{H}^{N-1} -\into \left(z, D \chi_{\{u > k\}}\right) \varphi + \lim_{n \to \infty} \into \left(z_n, D\chi_{\{ u_n > k\}}\right) \varphi .
		\end{multline}
		As an application of \eqref{stime norm l1 trac zn}, \eqref{dis per prov pair mis} and \eqref{desig}, we obtain
		\begin{equation*}
			\left|\int_{\partial \Omega}\varphi[z_k,\nu] \, \mathrm{d}\mathcal{H}^{N-1}  \right| \le 3\|f\|_{L^1(\Omega)} \|\varphi\|_{L^\infty(\partial \Omega)} \quad \text{for every $0 \le \varphi \in C^1(\Omega) \cap C(\overline{\Omega})$.}
		\end{equation*}
		We derive that
		\begin{equation*}
			\|[z_k,\nu]\|_{L^1(\partial \Omega)} \le 3 \|f\|_{L^1(\Omega)}.
		\end{equation*}
		Moreover, we observe that for almost every $h > k > 0$, the following property holds:
		\begin{equation*}
			\left|[z_k,\nu] \right|= \left|[z_k \chi_{\{u \le h\}},\nu]\right|\stackrel{\eqref{v fuori v dentro}}{=} \left|\chi_{\{u \le k\}}[z_h,\nu]\right| \le \left|[z_h,\nu]\right|,
		\end{equation*}
		because $z_k=u^m w \chi_{\{u \le k\}} $. This implies that $\left|[z_k,\nu]\right|$ is a nondecreasing sequence.
		Due to Beppo Levi's Theorem, there exists a function $\left\{z,\nu\right\} \in L^1(\partial \Omega)$ such that
		\begin{equation}\label{con trac z_k to z sul bordo}
			[z_k,\nu] \to \left\{z,\nu\right\} \quad \text{strongly in $L^1(\partial\Omega)$.}
		\end{equation}
		\medskip
		
		Now, we show that $\left\{z,\nu\right\}=[z,\nu]$ $\mathcal{H}^{N-1}$-a.e. on $\partial \Omega$. Let $\varphi \in C^1(\Omega)\cap C(\overline{\Omega})$. By the Gauss–Green formula \eqref{GG form}, we have
		\begin{equation}\label{per dim g la trac z}
			\into z_k \cdot \nabla \varphi -\int_{\partial \Omega} \varphi [z_k,\nu] \, \mathrm{d}\mathcal{H}^{N-1}= -\into \operatorname{div}z_k \varphi,
		\end{equation}
		we want to take limit as $k \to \infty$ in the previous equality.
		\\In the first integral on the left hand side of \eqref{per dim g la trac z}, we use Lebesgue's Theorem because $u$ is finite a.e. in $\Omega$ as stated in Remark \ref{rem u fqo}. For the second integral on the left hand side of \eqref{per dim g la trac z} we use \eqref{con trac z_k to z sul bordo}. On the right hand side, we note that
		\begin{equation*} 
		\begin{aligned}
				-\into \operatorname{div}z_k \varphi \stackrel{\eqref{ug pair con div lambda}}{=}& -\int_{\{u\le k\}} \operatorname{div}z \varphi - \into \left(z, D \chi_{\{ u \le k\}}\right) \varphi \\
			\stackrel{\eqref{sol 1}}{=} &\int_{\{u \le k\}} f \varphi + \into \left(z, D\chi_{\{u>k\}}\right).
		\end{aligned}
		\end{equation*}
		In particular, thanks to Lebesgue Theorem, we infer
		\begin{equation*} 
		\begin{aligned}
				-\lim_{k\to\infty}\into \operatorname{div} z_k \varphi
				= &\lim_{k \to \infty} \int_{\{u \le k\}} f \varphi + \into \left(z, D \chi_{\{u>k\}}\right) \varphi \\=&\into f \varphi,
		\end{aligned}
		\end{equation*}
		because \eqref{u fqo}  and \eqref{conv forte di zchiu>j} hold.
		\\Therefore, gathering together all the previous arguments, we obtain
		\begin{equation*}
			\into z \cdot \nabla \varphi - \int_{\partial \Omega} \varphi \left\{z,\nu\right\} \, \mathrm{d}\mathcal{H}^{N-1}=\into f\varphi.
		\end{equation*}
		As a consequence of \eqref{GG con z e reg}, it implies that
		\begin{equation*}
			\int_{\partial \Omega} \varphi \left\{z,\nu\right\} \mathrm{d}\mathcal{H}^{N-1}=\int_{\partial \Omega} \varphi [z,\nu] \mathrm{d}\mathcal{H}^{N-1} \quad \text{for every $\varphi \in C(\partial \Omega)$,}
		\end{equation*}
		which means that $\{z,\nu\}=[z,\nu]$ $\mathcal{H}^{N-1}$-a.e. on $\partial \Omega$. In particular, it holds
		\begin{equation}\label{conv trac z_k to z,nu}
			[z_k,\nu] \to [z,\nu] \quad \text{strongly in $L^1(\partial\Omega)$.}
		\end{equation}
		\medskip
		
		At this stage, we are in a position to prove the validity of the Gauss–Green formula \eqref{GG form definitiva}.
		Let us fix $0<a<b<\infty$ and $\alpha>0$, and let us take $k>b$. As a consequence of Lemma \ref{lem zab=z}, we have
		\begin{equation*}
			\begin{aligned}
				\into \left(z,DT_{a,b}(u)^\alpha\right) =& \lim_{k \to \infty} \into \left(z_k, DT_{a,b}(u)^\alpha\right) \\ \stackrel{\eqref{GG form}}{=} &\lim_{k \to \infty} \int_{\partial \Omega} T_{a,b}(u)^\alpha [z_k,\nu] \, \mathrm{d}\mathcal{H}^{N-1} - \into \operatorname{div}z_k \, T_{a,b}(u)^\alpha \\ = & \int_{\partial \Omega} T_{a,b}(u)^\alpha [z,\nu] \, \mathrm{d}\mathcal{H}^{N-1} - \into \operatorname{div}z \, T_{a,b}(u)^\alpha,
			\end{aligned}
		\end{equation*}
		where in the last equality we use \eqref{conv trac z_k to z,nu} and the identity
		$$\operatorname{div}z_k=\chi_{\{u\le k\}}\operatorname{div} z+(z,D\chi_{\{u\le k\}})$$
		jointly with \eqref{conv forte di zchiu>j}.
		This implies \eqref{GG form definitiva}.

		Finally, we prove \eqref{chi fuori chi dentro}. Let us take $\varphi \in C(\partial \Omega)$, $h>k>0$ and recalling that $z_k=z_k \chi_{\{u \le h\}}$, we get
		\begin{equation*}
			\begin{aligned}
				\int_{\partial \Omega} \varphi [z_k,\nu] \, \ensuremath{\mathrm d}\mathcal{H}^{N-1} = &\lim_{h \to \infty}\int_{\partial \Omega}\varphi [z_k \chi_{\{u \le h\}},\nu] \, \mathrm{d}\mathcal{H}^{N-1} \\ \stackrel{\eqref{v fuori v dentro}}{=}&\lim_{h \to \infty} \int_{\partial \Omega} \varphi \chi_{\{u \le k\}} [z_h,\nu] \, \mathrm{d}\mathcal{H}^{N-1}\\ =&\int_{\partial \Omega} \varphi \chi_{\{ u \le k\}}[z,\nu] \, \mathrm{d}\mathcal{H}^{N-1},
			\end{aligned}
		\end{equation*}
		where in the last equality we used \eqref{conv trac z_k to z,nu}. Since for every $\varphi \in C(\partial \Omega)$, it follows that $$\int_{\partial \Omega}\varphi [z_k,\nu] \, \mathrm{d}\mathcal{H}^{N-1}= \int_{\partial \Omega} \varphi \chi_{\{u \le k\}}[z,\nu]\, \mathrm{d}\mathcal{H}^{N-1},$$ and we deduce \eqref{chi fuori chi dentro}. This completes the proof.
	\end{proof}

	\begin{remark}\label{conv z_n a z}
		We explicitly stress that identities \eqref{chi fuori chi dentro} along with \eqref{Gallo} leads to
		\begin{equation*} 
			\lim_{n\to\infty}\int_\Omega \left(z_n, D\chi_{\{u_n>k\}}\right)\varphi=\int_\Omega \left(z, D\chi_{\{u>k\}}\right)\varphi,
		\end{equation*}
		for all $\varphi\in C^1(\Omega) \cap C(\overline\Omega)$. In particular, we get
		\[\lim_{n\to\infty}\int_\Omega \left(z_n, D\chi_{\{u_n>k\}}\right)=\int_\Omega \left(z, D\chi_{\{u>k\}}\right).\]
	\end{remark}
	
	The next result is an inequality that proves be useful in establishing the boundary condition.
	\begin{lemma}\label{casi frontera}
		Assume $m>0$ and $0 \le f \in L^1(\Omega)$. Consider $u_n$ solution to problem \eqref{Prob trans tronc} with associated vector field $w_n$. Then, for every $q>0$ it holds
		\begin{equation}\label{dis per prov sol 3}
			\frac{T_{a,b}(u)^{m(q+1)} }{q+1} \le - \frac{T_{a,b}(u)^{mq} }{q}[z,\nu] \quad \text{for $\mathcal{H}^{N-1}$-a.e. $x \in \partial \Omega$,}
		\end{equation}
		for a.e. $0<a<b<\infty$.
	\end{lemma}
	\begin{proof}
		It is an adaptation of \cite[Lemma 5.7]{BOPS}, however, for the sake of completeness, we outline the key details.
		\\Let us fix $q>0$ and let us take $\frac{T_{a,b}(u_n)^{mq}}{q} \varphi \in BV(\Omega)\cap L^\infty(\Omega)$ with $0 \le \varphi \in C^1(\overline{\Omega})$  as a test function in \eqref{exp f test n}, obtaining
		\begin{equation}\label{per il bordo 1}
			\begin{aligned}
				&\into \frac{1}{q}\left(z_n, DT_{a,b}(u_n)^{mq} \right) \varphi + \frac{1}{q}\into T_{a,b}(u_n)^{mq} z_n \cdot \nabla \varphi \\- &\int_{\partial \Omega} \frac{T_{a,b}(u_n)^{mq}}{q} \varphi [z_n, \nu] \, \ensuremath{\mathrm d}\mathcal{H}^{N-1}= \into f_n \frac{T_{a,b}(u_n)^{mq}}{q} \varphi.
			\end{aligned}
		\end{equation}
		We want to take the limit as $n$ tends to infinity in this equation, but before we will perform some manipulations.
		For the first integral on the left hand side of \eqref{per il bordo 1}, we note that
		\begin{equation}\label{andando al bordo 1}
			\begin{aligned}
				\into \frac{1}{q}\left(z_n, DT_{a,b}(u_n)^{mq}\right) \varphi \stackrel{\eqref{pair con der teta}}{=} &\into \frac{1}{q} \theta\left(z_n, DT_{a,b}(u_n)^{mq},x\right) \left|DT_{a,b}(u_n)^{mq} \right| \varphi \\ \stackrel{\eqref{eq 4.5 cdc}}{=} &\into \frac{1}{q}\theta\left(z_n, DT_{a,b}(u_n),x\right) |DT_{a,b}(u_n)^{mq}| \varphi \\ \stackrel{\eqref{ch wj}}{=} &\into \theta \left(z_n, DT_{a,b}(u_n),x\right) mT_{a,b}(u_n)^{mq-1}|DT_{a,b}(u_n)| \varphi \\ \stackrel{\eqref{pair con der teta}}{=} &\into mT_{a,b}(u_n)^{mq-1} \left(z_n, DT_{a,b}(u_n)\right) \varphi \\ \stackrel{\eqref{sol 2 n}}{=} &\frac{m}{m+1}\into T_{a,b}(u_n)^{mq-1} \left|DT_{a,b}(u_n)^{m+1}\right| \varphi \\ \stackrel{\eqref{ch wj}}{=} &\frac{1}{q+1}\into \left|DT_{a,b}(u_n)^{m(q+1)} \right| \varphi.
			\end{aligned}
		\end{equation}
		Regarding the third integral on the left hand side of \eqref{per il bordo 1}, we employ $T_{a,b}(u_n)^m \le - [z_n, \nu]$ for a.e. $0<a<b<\infty$ and for $\mathcal{H}^{N-1}$-a.e. $x \in \partial \Omega$ (by \cite[Lemma 5.8]{BOPS}) and deduce that
		\begin{equation}\label{andando al bordo 2}
			\begin{aligned}
				-&\int_{\partial \Omega} \frac{1}{q}T_{a,b}(u_n)^{mq} \varphi [z_n, \nu]
				\, \ensuremath{\mathrm d} \mathcal{H}^{N-1}\\
				\ge &\int_{\partial \Omega} \frac{1}{q}T_{a,b}(u_n)^{mq} T_{a,b}(u_n)^m \varphi \, \ensuremath{\mathrm d} \mathcal{H}^{N-1}  \\
				\ge &\int_{\partial \Omega} \frac{1}{q+1} T_{a,b}(u_n)^{m(q+1)} \varphi \, \ensuremath{\mathrm d}\mathcal{H}^{N-1}.
			\end{aligned}
		\end{equation}
		Putting together \eqref{andando al bordo 1} and \eqref{andando al bordo 2} in \eqref{per il bordo 1}, it yields
		\begin{equation}\label{per il bordo 2}
			\begin{aligned}
				\frac{1}{q+1}&\into \left|D T_{a,b}(u_n)^{m(q+1)}\right| \varphi + \frac{1}{q} \into T_{a,b}(u_n)^{mq} z_n \cdot \nabla \varphi \\ + \frac{1}{q+1} &\int_{\partial \Omega} T_{a,b}(u_n)^{m(q+1)}  \varphi \, \ensuremath{\mathrm d}\mathcal{H}^{N-1} \le \into f_n \frac{T_{a,b}(u_n)^{mq}}{q} \varphi.
			\end{aligned}
		\end{equation}
		Next we take the limit as $n$ tends to infinity. For the first and the third integral on the left hand side of \eqref{per il bordo 2} we use \eqref{sci grad con bordo} recalling that \eqref{conv in L^1} holds. On the other hand, for the second integral on the left hand side of \eqref{per il bordo 2} we use properties \eqref{con deb zn to z} and \eqref{conv TAB}. Finally for the integral on the right hand side we use Lebesgue's Theorem. Hence we get
		\begin{equation*}
			\begin{aligned}
				\frac{1}{q+1}&\into \left|D T_{a,b}(u)^{m(q+1)}\right| \varphi + \frac{1}{q} \into T_{a,b}(u)^{mq} z \cdot \nabla \varphi \\+ \frac{1}{q+1} &\int_{\partial \Omega} T_{a,b}(u)^{m(q+1)}  \varphi \, \ensuremath{\mathrm d}\mathcal{H}^{N-1} \le \into f \frac{T_{a,b}(u)^{mq}}{q} \varphi.
			\end{aligned}
		\end{equation*}
		Using \eqref{GG form definitiva}, we obtain
		\begin{equation}\label{per il bordo 3}
			\begin{aligned}
				\frac{1}{q+1}&\into \left|D T_{a,b}(u)^{m(q+1)}\right| \varphi + \frac{1}{q+1} \int_{\partial \Omega} T_{a,b}(u)^{m(q+1)}  \varphi \, \ensuremath{\mathrm d}\mathcal{H}^{N-1} \\\le \frac{1}{q}&\into \left(z, DT_{a,b}(u)^{mq}\right) \varphi - \frac{1}{q}\int_{\partial \Omega} T_{a,b}(u)^{mq} \varphi [z,\nu] \, \ensuremath{\mathrm d}\mathcal{H}^{N-1}.
			\end{aligned}
		\end{equation}
		Recalling Lemma \ref{lem zab=z} and repeating the same arguments of \eqref{andando al bordo 1}, we get
		\begin{equation*}
			\frac{1}{q}\into \left(z, DT_{a,b}(u)^{mq}\right) \varphi = \frac{1}{q+1} \into \left|DT_{a,b}(u)^{m(q+1)}\right| \varphi.
		\end{equation*}
		Therefore \eqref{per il bordo 3} becomes
		\begin{equation*}
			\frac{1}{q+1}\int_{\partial \Omega} T_{a,b}(u)^{m(q+1)} \varphi \, \ensuremath{\mathrm d}\mathcal{H}^{N-1} \le  - \frac{1}{q}\int_{\partial \Omega} T_{a,b}(u)^{mq} \varphi [z,\nu] \, \ensuremath{\mathrm d}\mathcal{H}^{N-1},
		\end{equation*}
		which implies \eqref{dis per prov sol 3} and this concludes the proof.
	\end{proof}
	
	The boundary condition \eqref{sol 3} now easily follows from Lemma \ref{casi frontera}.
	\begin{lemma}\label{lem cond bord}
		Assume $m>0$ and $0 \le f \in L^1(\Omega)$. Let us consider $u_n$ a solution to problem \eqref{Prob trans tronc} with associated vector field $w_n$. Then, it holds \eqref{sol 3}.
	\end{lemma}
	\begin{proof}
		It is an adaptation of \cite[Lemma 5.8]{BOPS}, but for the sake of completeness we underline the main steps.
		We fix $0<a<b<\infty$ such that \eqref{prob con t e fi} holds.
		\\Dividing both sides of \eqref{dis per prov sol 3} by $T_{a,b}(u)^{mq}$, we obtain
		\begin{equation*}
			\frac{q}{q+1} T_{a,b}(u)^m \le - [z,\nu] \quad \text{for $\mathcal{H}^{N-1}$-a.e. $x \in \partial \Omega$,}
		\end{equation*}
		hence, letting $q$ go to infinity, it yields
		\begin{equation*} 
			T_{a,b}(u)^m \le - [z, \nu] \quad \text{for $\mathcal{H}^{N-1}$-a.e. $x \in \partial \Omega$.}
		\end{equation*}
		Since it holds for almost all $0<a<b<\infty$,
		taking the limit as $a$ tends to $0$ and $b$ tends to infinity, we get
		\begin{equation}\label{1 per sol 3}
			\left(u^\Omega\right)^m \le -[z, \nu] \quad \text{for $\mathcal{H}^{N-1}$-a.e. $x \in \partial \Omega \cap \{u>0\}$.}
		\end{equation}
		This implies that $\left(u^\Omega\right)^m \in L^1(\partial \Omega)$, because $\left(u^\Omega\right)^m$ is nonnegative.
		\\Let us show the reverse inequality. To this end, fix again $0<a<b<\infty$, recalling that $z_{a,b}, w \chi_{\{a <u<b\}} \in \DM(\Omega)$, it follows that
		\begin{equation*}
			\left|[z,\nu] \chi_{\{a<u<b\}}\right| \stackrel{\eqref{chi fuori chi dentro}}{=} \left|[z_{a,b},\nu] \right| \stackrel{\eqref{v fuori v dentro}}{\le} T_{a,b}(u)^m,
		\end{equation*}
		Letting $a$ go to 0 and $b$ to infinity in $\partial \Omega \cap \{u>0\}$, we have
		\begin{equation}\label{2 per sol 3}
			|[z,\nu]| \le \left(u^\Omega\right)^m \quad \text{for $\mathcal{H}^{N-1}$-a.e. $x \in \partial \Omega \cap \{u>0\}$.}
		\end{equation}
		Putting together \eqref{1 per sol 3} and \eqref{2 per sol 3}, we obtain \eqref{sol 3} and this concludes the proof.
	\end{proof}
	We are now in a position to prove our main result.
	\begin{proof}[End of proof of Theorem \ref{teo tras 1}]
		The proof follows as a direct consequence of the preceding results. In Corollary \ref{cor ext u}, we establish the existence of a Lebesgue measurable, nonnegative function $u$ such that, for almost every $0 < a < b < \infty$, the truncated function $T_{a,b}(u)$ belongs to $BV(\Omega)$. In particular, as noted in Remark \ref{rem u fqo}, we get that $u(x) < \infty$ for almost every $x \in \Omega$.
		
		Moreover, Lemma \ref{lem u j v} ensures that $u \in DTBV(\Omega)$, providing further insight into its regularity. In addition, Lemma \ref{lem ext z} guarantees the existence of a vector field $w \in L^\infty(\Omega)^N$ with $\|w\|_{L^\infty(\Omega)^N} \leq 1$, from which we define $z := u^m w \in L^1(\Omega)^N$. This construction yields the fundamental identity \eqref{sol 1}, implying that $z \in \mathcal{DM}^1(\Omega)$.
		
		Furthermore, Lemma \ref{lem per sol 2} establishes the validity of \eqref{sol 2}. Finally, Lemma \ref{lem trac z in l1} ensures that $[z, \nu] \in L^1(\partial \Omega)$, while Lemma \ref{lem cond bord} confirms the boundary condition \eqref{sol 3}, thereby completing the proof.
	\end{proof}
	
	\section{Regularity results}\label{sec6}
	We finish this work by studying the summability of the solution that occurs when the datum is between the Lebesgue spaces $L^1(\Omega)$ (our setting) and $L^N(\Omega)$ (the framework of \cite{BOPS}). This regularity result is consistent with the findings presented.
\begin{theorem}
	Assume $m>0$ and $0 \le f \in L^p(\Omega)$, with $1<p<N$. Then the solution $u$ to problem \eqref{Prob trans tronc} we have found in Theorem \ref{teo tras 1} satisfies
	\begin{equation}\label{reg 1 ris}
		u^m \in L^{\frac{Np}{N-p}}(\Omega).
	\end{equation}
	Moreover, if $p \ge \frac{N(m+1)}{Nm+1}$, then it holds
	\begin{equation}\label{reg 2 ris}
		u^{m+1} \in BV(\Omega).
	\end{equation}
\end{theorem}
\begin{proof}
	Let $u_n$ be an approximate solution with associated vector field $w_n$.
	We choose $T_{a,\infty}(u_n)^\alpha$, with $a > 0$, as a test function in problem \eqref{Prob trans tronc}, $\alpha > 0$ to be determined. This yields
	$$
	\int_\Omega \left(z_n, D T_{a,\infty}(u_n)^\alpha\right) - \int_{\partial \Omega} T_{a,\infty}(u_n)^\alpha [z_n,\nu] \, \mathrm{d}\mathcal{H}^{N-1} = \int_\Omega f_n T_{a,\infty}(u_n)^\alpha.
	$$
	
	Applying \cite[Lemma 2.11]{BOPS} and H\"older's inequality, we obtain
	\begin{equation}\label{per la seconda reg}
		\frac{\alpha}{\alpha+m} \|T_{a,\infty}(u_n)^{\alpha + m}\|_{BV(\Omega)} \le \int_\Omega f_n T_{a,\infty}(u_n)^\alpha
		\le \|f\|_{L^p(\Omega)} \|T_{a,\infty}(u_n)^\alpha\|_{L^{\frac{p}{p-1}}(\Omega)},
	\end{equation}
	where we recall that $T_{a,\infty}(u_n)^m \le -[z_n,\nu]$, as shown in \cite[Lemma 5.8]{BOPS}.
	
	\medskip
	
	Next, Sobolev's embedding \eqref{sob emb} implies
	\begin{equation}\label{per reg 1}
		\frac{\alpha}{\alpha + m} \mathcal{S}_1^{-1} \|T_{a,\infty}(u_n)^{\alpha+m}\|_{L^{\frac{N}{N-1}}(\Omega)} \le \|f\|_{L^p(\Omega)} \|T_{a,\infty}(u_n)^\alpha\|_{L^{\frac{p}{p-1}}(\Omega)}.
	\end{equation}
	We now select $\alpha > 0$ such that
	\begin{equation}\label{la scelta di alpha}
		(\alpha + m) \frac{N}{N-1} = \alpha \frac{p}{p - 1} \quad \Rightarrow \quad \alpha = \frac{(p - 1)mN}{N - p}.
	\end{equation}
	Substituting this choice into \eqref{per reg 1}, we deduce
	
	\begin{equation}\label{penultimo}
		\left( \int_\Omega T_{a,\infty}(u_n)^{\frac{mNp}{N - p}} \right)^{\frac{N - p}{Np}} \le \frac{p(N - 1)}{(p - 1)N} \mathcal{S}_1 \|f\|_{L^p(\Omega)}.
	\end{equation}
	
	As a consequence of Corollary \ref{cor ext u}, we can pass to the limit first as $n \to \infty$ and then as $a \to 0$, thus obtaining
	$$
	\|u^m\|_{L^{\frac{Np}{N - p}}(\Omega)} \le \frac{p(N - 1)}{(p - 1)N} \mathcal{S}_1 \|f\|_{L^p(\Omega)},
	$$
	which establishes the validity of \eqref{reg 1 ris}.
	
	\medskip
	
	To finish the proof, we now verify property \eqref{reg 2 ris}. Having in mind the usual embeddings satisfy by Lebesgue spaces, we may assume $p = \frac{N(m + 1)}{Nm + 1}$. This choice of $p$ implies $\alpha =1$ in \eqref{la scelta di alpha}. We deduce from  \eqref{per la seconda reg} that
	\begin{equation}\label{ultimo}
		\|T_{a,\infty}(u_n)^{m+1}\|_{BV(\Omega)} \le \frac{p(N - 1)}{(p - 1)N} \|f\|_{L^p(\Omega)} \|T_{a,\infty}(u_n)\|_{L^{\frac{p}{p-1}}(\Omega)}.
	\end{equation}
	Since $\frac{p}{p-1}=\frac{Npm}{N-p}$, it follows from \eqref{penultimo} that the right hand side of \eqref{ultimo} is bounded. The lower semicontinuity of the $BV$-norm allows us to take the limit as $n$ goes to infinity, it results
	$$
	\|T_{a,\infty}(u)^{m+1}\|_{BV(\Omega)} \le \mathcal{C},
	$$
	where $\mathcal{C}>0$ is a constant independent of $a$. Hence $u^{m+1}\in BV(\Omega)$ and we are done.
	\end{proof}
	
	\begin{remark}
		We point out that, letting $p$ go to 1 in \eqref{reg 1 ris}, we obtain the Lebesgue space $L^{\frac{N}{N-1}}(\Omega)$ instead of the Marcinkiewicz space $L^{\frac{N}{N-1},\infty}(\Omega)$ (see Remark \ref{rem u fqo}).  In this way, we ``almost'' get the actual regularity for $u^m$. This is a usual feature when studying elliptic equations with $L^1$ (or measure) data: the expected Lebesgue space is not attained and instead the solution lies in the corresponding Marcinkiewicz space (see, for instance, \cite{S} for the linear setting or \cite{BG} in a nonlinear case). In the other extreme case, we have that when $p$ tends to $N$, we arrive at $u^m\in L^\infty(\Omega)$, which is the summability proven in \cite{BOPS}.
		
		On the other hand, the parameter $p=\frac{N(m + 1)}{Nm + 1}$ is exactly that found in \cite{BOPS} to have a $BV$-estimate (see \cite[Lemma 4.2]{BOPS}).
	\end{remark}

	\section*{Acknowledgments}
	F. Balducci is partially supported by the Gruppo Nazionale per l’Analisi Matematica, la Probabilità e le loro Applicazioni (GNAMPA) of the Istituto Nazionale di Alta Matematica (INdAM).
	S. Segura de Le\'on has partially been supported by Grant PID2022-136589NB-I00 funded by MCIN/
	AEI/10.13039/501100011033 and by Grant RED2022-134784-T also funded by MCIN/AEI/10.13039/
	501100011033.

\end{document}